%% file: author_arxiv.tex
\documentclass[a4paper,10pt]{amsart}
\usepackage[margin=3cm]{geometry}
\usepackage[foot]{amsaddr}
\usepackage{amsmath,epsfig,multirow}
\usepackage{mathtools}
\usepackage{amsmath, bm}
\usepackage{bm}
\usepackage{pgfplots}
\usepackage{mathrsfs}
\usepackage{systeme}
\usepackage{amsfonts}
\usepackage{subfigure}
\usepackage{graphicx,subfigure}
\usepackage{hyperref}
\usepackage{mathrsfs}

\usepackage{graphicx}        
\usepackage{multicol}        
\usepackage[bottom]{footmisc}

\graphicspath{./figures}
\makeindex             
\begin{document}

\title[Non-intrusive PCE in CFD Problems and Comparison to POD]{Non-Intrusive Polynomial Chaos Method Applied to Full-Order and Reduced Problems in
Computational Fluid Dynamics: a Comparison and Perspectives}

\author{Saddam Hijazi$^1$}
\author{Giovanni Stabile$^1$}
\author{Andrea Mola$^1$}
\author{Gianluigi Rozza$^1$}

\address{$^1$ mathLab, Mathematics Area, SISSA, via Bonomea 265, I-34136 Trieste, Italy}
%
%
\maketitle

\begin{abstract}
{\color{black}In this work, Uncertainty Quantification (UQ) based on non-intrusive Polynomial Chaos Expansion (PCE) is applied to the CFD problem of the flow past an airfoil with parameterized angle of attack and inflow velocity. To limit the computational cost associated with each of the simulations required by the non-intrusive UQ algorithm used, we resort to a Reduced Order Model (ROM) based on Proper Orthogonal Decomposition (POD)-Galerkin approach. A first set of results is presented to characterize the accuracy of the POD-Galerkin ROM developed approach with respect to the Full Order Model (FOM) solver (OpenFOAM). A further analysis is then presented to assess how the UQ results are affected by substituting the FOM predictions with the surrogate ROM ones.}
\end{abstract}

\section{Introduction}
\label{sec:intro}
Many methods have been developed to assess how uncertainties of input parameters propagate, through Computational Fluid Dynamics (CFD) numerical simulations, into the outputs of interest. The aim of this work is to carry out a study on the application of non-intrusive Polynomial Chaos Expansion (PCE) to CFD problems. {\color{black}The PCE method is a way of representing random variables or random processes in terms of orthogonal polynomials. One important feature of PCE is the possibility of decomposing the random variable into separable deterministic and stochastic components \cite{loeve1978probability,hosder2006non}}. By a computational stand point, the main problem in PCE consists in finding the deterministic coefficients of the expansion. In  non-intrusive PCE, no changes are made in the simulations code, and the coefficients are computed in a post processing phase which follows the simulations. Thus, the deterministic terms in the expansion are obtained via a sampling based approach such as the one used in \cite{isukapalli1999uncertainty,reagana2003uncertainty}. In this framework, samples of the input parameters are prescribed and then numerical simulations are carried out for each sample.

{\color{black}Once the output of the simulations corresponding to each sample is evaluated, it is used to obtain the PCE coefficients. In the projection approach, the orthogonality of the polynomials is exploited to compute the deterministic coefficients in the expansion through integrals in the sampling space. As the sampling points chosen are quadrature points for such integrals, the computational cost will grow exponentially as the parameter space dimension increases. This is of course quite undesirable, given the considerable computational cost of the CFD simulations associated to the output evaluation at each sampling point. To avoid such problem, in this work the PCE expansion coefficients are computed using a regression approach which is based on least squares minimization}.\par
{\color{black} To explore even further reductions of the computational cost associated with sample points output evaluations, in the present work we apply the PCE algorithm both to the full order CFD model and to a reduced order model based on POD-Galerkin approach. In the last decade, there have been several efforts to develop reduced order models and apply them to industrial continuous mechanics problems governed by parameterized PDEs. We refer the interested readers to \cite{hesthaven2015certified,quarteroniRB2016,chinesta2017model} for detailed theory on ROMs for parameterized PDE problems. More in particular, the solution of parameterized Navier--Stokes problem in a reduced order setting is discussed in \cite{Quarteroni2007}. In such work, the FOM discretization was based on Finite Element Method (FEM). Projection-based Reduced Order Methods (ROMs) have in fact been mainly developed for FEM, but in the last years many efforts have been dedicated to extend them to Finite Volume Method (FVM) and to CFD problems with high Reynolds numbers. Some examples of the application of ROMs based on the Reduced Basis (RB) method to a finite volume setting are found in \cite{Haasdonk2008277,HOR08,Drohmann2012}. In this work, we instead focus on POD-Galerkin methods applied to CFD computations based on FVM discretization. A large variety of works related to POD-Galerkin can be found in the literature, and here we refer only to some of them \cite{noack1994,Akhtar2009,Bergmann2009516,Kunisch2002492,Burkardt2006337,Baiges2014189}, as examples. As for POD-Galerkin approach applied to Navier--Stokes flows discretized via FVM, we mention \cite{Lorenzi2016151}, in which the authors treat the velocity pressure coupling in the reduced model using the  same set of coefficients for both velocity and pressure fields. In \cite{Stabile2017}, the coefficients of velocity and pressure are instead different, and Poisson equation for pressure is added to close the system at reduced order level. In \cite{STABILE2018} a stabilization method for the finite volume ROM model is presented. In \cite{Carlberg2018} a study on conservative reduced order model for finite volume method is discussed. For applications of ROMs to UQ problems we refer the readers to \cite{chen2017reduced,chen2014comparison,gunzburger2007reduced}.} \par
 PCE is a tool that is independent of the output evaluator and in this work we will apply it {\color{black} to output parameters both obtained from the full order solution and to its POD-Galerkin reduced order counterpart}. In this regard, the objective of the present work is to assess whether PCE {\color{black} results are significantly} influenced by the use of a POD-Galerkin based model reduction approach. To this end, we will apply POD model reduction to CFD simulations based on incompressible steady Navier--Stokes equations, and compare the PCE coefficients and sensitivities obtained for the reduced order solution to the ones resulting from the full order simulations.\par
{\color{black} 
This article is organized as follows: in section \ref{sec:HF_model} the physical problem under study is described at the full order level. In section \ref{sec:ROM} the  reduced order model is introduced. In particular, the most relevant notions on projection based methods are reported in subsection \ref{sec:project}, while boundary conditions treatment is discussed in subsection \ref{sec:boundary}. The theory of the non intrusive PCE is summarized in section \ref{sec:PCE} with direct reference to the quantities of interest in the present work. Numerical results are presented in section \ref{sec:results}, starting with the ones of the reduced order model in subsection \ref{sec:ROM_res} and then the PCE results in \ref{sec:PCE_res}. Finally, conclusions and possible directions of future work are discussed in section \ref{sec:concl}. }

\section{The Physical Problem}
\label{sec:HF_model}
\begin{figure}
\centering
\includegraphics[width=0.7\textwidth]{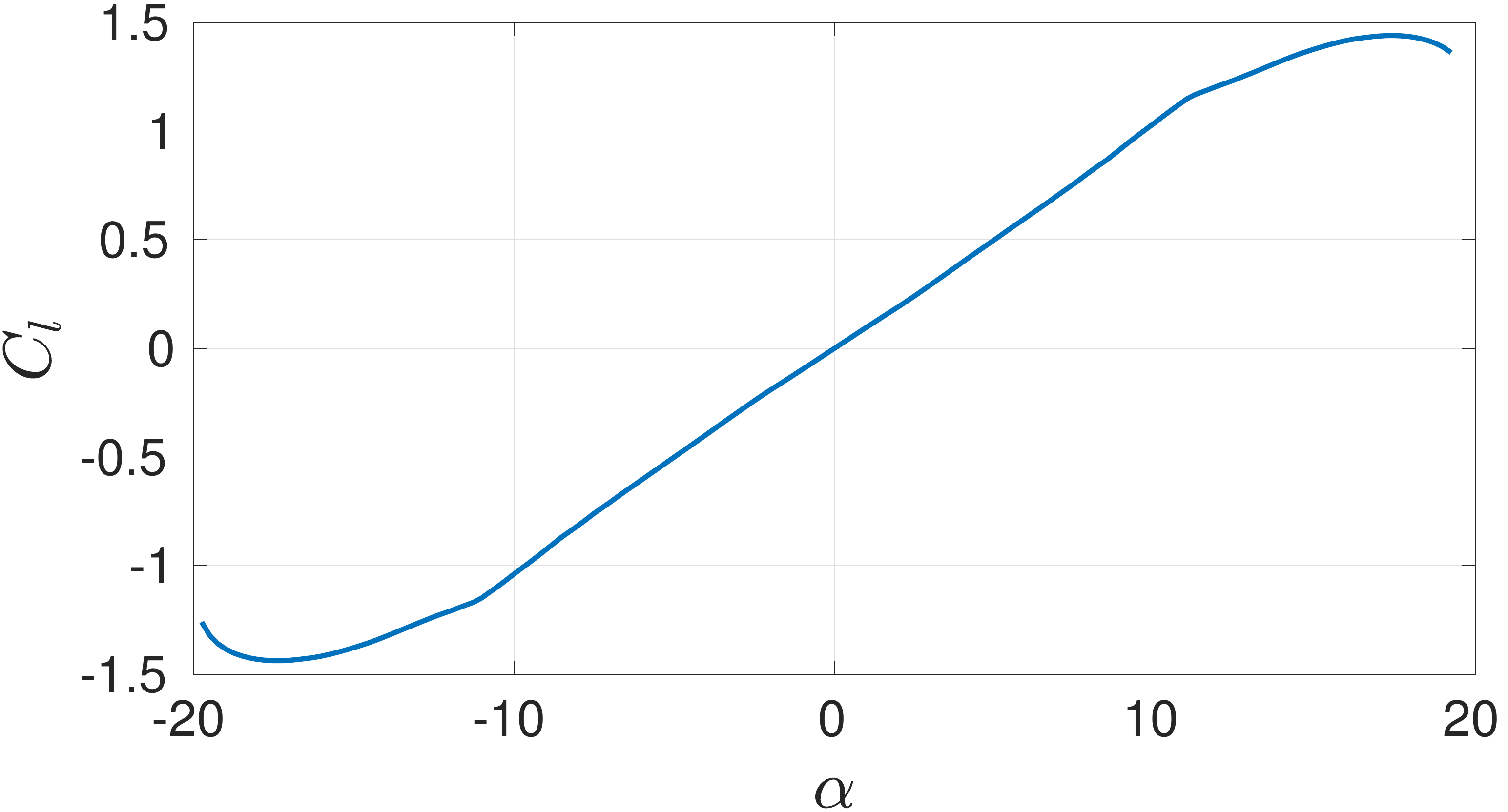}
\caption{The lift coefficient curve for the airfoil NACA0012.}\label{fig:cl_curve}
\end{figure}

In this section, we describe the physical problem of interest which consists into the flow around an airfoil subjected to variations of the angle of attack {\color{black} and inflow velocity}. In aerospace engineering, the angle of attack is the angle that lies between the flow velocity vector at infinite distance from the airfoil ($U_\infty$) and the chord of the airfoil, see \autoref{fig:aoa}. We are interested in finding the angle of attack that produces the maximum lift coefficient before stall happens. {\color{black} \autoref{fig:cl_curve} depicts for the lift coefficient curve of the airfoil NACA $0012-64$ \cite{abbott1999theory, Sarkar2009} at a fixed Reynolds number of $10^6$. The plot suggests that as the angle of attack increases, the  lift coefficient grows until flow separation occurs leading to a loss of the lift force. At laminar flow regimes such as the ones that will be analyzed later in section \ref{sec:results}, such stall phenomenon happens in a mild fashion, as opposed to more abrupt stalls observed at higher velocities like the one in \autoref{fig:cl_curve}. It can be noticed from the plot that the lift coefficient reaches its maximum value when the angle is about $17$ degrees before stall happens. For lower Reynolds numbers, the maximum in the $C_L$-$\alpha$ curve is observed at higher angles of attack.
\begin{figure}
\centering
\def\svgwidth{0.5\linewidth}
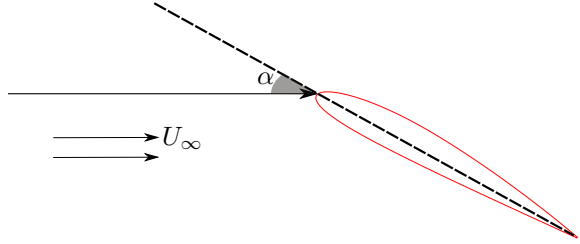
\caption{The angle of attack on an airfoil.}\label{fig:aoa}
\end{figure}
The fluid dynamic problem is mathematically governed by the steady Navier-Stokes equations which read as follows:} 
\begin{equation}\label{eq:navstokes}
\begin{cases}
(\bm{u} \cdot \bm{\nabla}) \bm{u}- \bm{\nabla} \cdot \nu \bm{\nabla} \bm{u}=-\bm{\nabla}p &\mbox{ in } \Omega_f, \\
\bm{\nabla} \cdot \bm{u}=\bm{0} &\mbox{ in } \Omega_f, \\
\bm{u} (x) = \bm{f}(\bm{x},\bm{\mu}) &\mbox{ on } \Gamma_{In}, \\
\bm{u} (x) = \bm{0} &\mbox{ on } \Gamma_{0}, \\ 
(\nu\nabla \bm{u} - p\bm{I})\bm{n} = \bm{0} &\mbox{ on } \Gamma_{Out}, \\
\end{cases}
\end{equation}
where $\Gamma = \Gamma_{In} \cup \Gamma_0 \cup \Gamma_{Out}$ is the boundary of the fluid domain $\Omega_f$ and is composed by three different parts $\Gamma_{In}$, $\Gamma_{Out}$ and $\Gamma_0$, indicating respectively inlet boundary, outlet boundary and wing. In the flow equations $\bm{u}$ is the flow velocity vector, $\nu$ is the fluid kinematic viscosity, and $p$ is the normalized pressure, which is divided by the fluid density $\rho_f$. As for the boundary conditions in which $\bm{f}$ is a generic function that prescribes the value of the velocity on the inlet $\Gamma_{In}$ and it is parameterized through the vector quantity $\bm{\mu}$. In the present work the problem is solved using a finite volume discretization technique \cite{Versteeg1995AnIT, Moukalled:2015:FVM:2876154,jasak1996error,ohlbergerFV}, the standard approach is to work with a Poisson equation for pressure, rather than directly with the continuity equation. System \eqref{eq:navstokes} is then modified into:
{\color{black}
\begin{equation}\label{eq:navstokesFV}
\begin{cases}
(\bm{u} \cdot \bm{\nabla}) \bm{u}- \bm{\nabla} \cdot \nu \bm{\nabla} \bm{u}=-\bm{\nabla}p &\mbox{ in } \Omega_f, \\
\Delta p = - \nabla \cdot (\bm{u} \cdot \nabla ) \bm{u} &\mbox{ in } \Omega_f, \\
\bm{u} (x) = \bm{f}(x,\mu) &\mbox{ on } \Gamma_{In}, \\
\bm{u} (x) = \bm{0} &\mbox{ on } \Gamma_{0}, \\ 
(\nabla \bm{u}) \bm{n} = \bm{0} &\mbox{ on } \Gamma_{Out}, \\ 
\nabla p \cdot \bm{n} = 0 &\mbox{ on } \Gamma \setminus \Gamma_{Out}, \\ 
p = 0 &\mbox{ on } \Gamma_{Out} .          
\end{cases}
\end{equation}}
In the above system of equations all the quantities assume the same meaning of those presented in \eqref{eq:navstokes}. {\color{black}The Poisson equation for pressure is obtained taking the divergence of the momentum equation, and then exploiting the divergence free constraint on velocity. The two equations are solved in a segregated fashion, making use of the SIMPLE algorithm \cite{PATANKAR19721787}.} Historically, the FVM discretization technique has been widely used in industrial applications and for flows characterized by higher values of the \emph{Reynolds} number. One important feature of the FVM is that it ensures that conservative laws are satisfied at local level. In this work the Full Order Model (FOM) simulations are carried out making use of the finite volume open source C\texttt{++} library OpenFOAM\textsuperscript{\textregistered} (OF) \cite{weller1998tensorial}. 


\section{The Reduced Order Model}\label{sec:ROM}
{\color{black} The FOM simulations carried out by OpenFOAM present a high computational cost. In the framework of a many query problem such as the one associated with non-intrusive PCE employed in this work, the search for ways to reduce the computational cost becomes paramount. For this reason, we resort to reduce order modelling and we couple it with PCE in the next sections. In This section, we recall the notion of ROM and the POD approach to build the reduced order spaces. Here, only few details are addressed, while for further information on how to adapt ROM for finite volume discretization method the reader may refer to \cite{Lorenzi2016151,Stabile2017,STABILE2018}.


The key assumption of ROMs is that one can find a low dimensional space in which it is possible to express the solution of the full order problem with good approximation properties. That space is spanned by the reduced order modes \cite{hesthaven2015certified}. The latter assumption translates to the following decomposition of the velocity and pressure fields:
\begin{eqnarray}\label{eq:vel_exp}
\bm{u}(\bm{x},\bm{\mu}) \approx \bm{u}_r(\bm{x},\bm{\mu}) = \sum_{i=1}^{N_u} a_i(\bm{\mu}) \bm{\phi}_i(\bm{x}),
\\ \label{eq:pre_exp}
p(\bm{x},\bm{\mu})\approx p_r(\bm{x},\bm{\mu}) = \sum_{i=1}^{N_p} b_i (\bm{\mu}) {\chi_i}(\bm{x}) ,
\end{eqnarray}
where $\bm{u}_r(\bm{x},\bm{\mu})$ and $p_r(\bm{x},\bm{\mu})$ are the reduced order approximations of velocity and pressure, respectively, $a_i$ and $b_i$ are scalar coefficients that depend on the parameter value $\bm{\mu}$, $\bm{\phi}_i$ and ${\chi_i}$ are the basis functions of the reduced basis spaces for velocity and pressure, respectively. $N_u$ and $N_p$ represent the dimension of the reduced basis spaces for velocity and pressure, respectively, obviously $N_u$ and $N_p$ are not supposed to have the same value.\par 
The next step in constructing the reduced order model is to generate the reduced order space. For such step we resort to a POD approach. The POD space is constructed by solving the following minimization problem:
\begin{equation}\label{eq:min1}
\mathbb{V}_{POD}= \mbox{arg min} \frac{1}{N_{s}} \sum_{n=1}^{N_s} ||  \bm{u}_n- \sum_{n=1}^{N_s}(\bm{u}_n, \bm{\phi}_i )_{L^2(\Omega)} \bm{\phi}_i ||_{L^2(\Omega)}^2,
\end{equation}
where $\bm{u}_n$ is a solution snapshot obtained for a certain parameter value $\bm{\mu_n}$ and $N_s$ is the total number of solution snapshots. One can see that the reduced order space (or $\mathbb{V}_{POD}$) is optimal in the sense that it is spanned by the modes that minimize the projection error between the fields and their projection into the modes. For further details on how the problem \ref{eq:min1} is solved one can refer to \cite{Stabile2017}.}

\subsection{Projection based ROM}\label{sec:project}
The next step in building the reduced order model (this procedure is referred as POD-Galerkin projection) is to project the momentum equation of \eqref{eq:navstokes} onto the POD space spanned by the POD velocity modes, namely:
\begin{equation}\label{eq:l2proj_vel}
\left( \bm{\phi}_i,(\bm{u} \cdot \nabla ) \bm{u} -\nu \Delta \bm{u} +\nabla p \right)_{L^2(\Omega)} = 0 .
\end{equation}
Inserting the approximations \eqref{eq:vel_exp} and \eqref{eq:pre_exp} into \eqref{eq:l2proj_vel} yields the following reduced system:
\begin{equation}\label{eq:non_linear_sys_lam}
\nu \bm{B} \bm{a} - \bm{a}^T \bm{C} \bm{a}-\bm{H}\bm{b} = \bm{0},
\end{equation}
where $\bm{a}$ and $\bm{b}$ are the vectors of coefficients for reduced velocity and reduced pressure respectively, while $\bm{B},\bm{C},\bm{H}$ are the reduced discretized differential operators which are computed as follows :
\begin{align}
& B_{ij}=\left( \bm{\phi}_i ,\Delta \bm{\phi}_j\right)_{L^2(\Omega)}, \\
& C_{ijk}=\left( \bm{\phi}_i ,(\bm{\phi}_j \cdot \bm{\nabla}) \bm{\phi}_k)\right)_{L^2(\Omega)} , \label{eq:div_phi}\\
& H_{ij} = \left( \bm{\phi}_i , \nabla \chi_j \right)_{L^2(\Omega)}.
\end{align}

To solve the system \eqref{eq:non_linear_sys_lam}, one needs $N_p$ additional equations. The continuity equation cannot be directly used because the snapshots are divergence free and so are the velocity POD modes. The available approaches to tackle this problem are either the use of the Poisson equation \cite{Stabile2017,STABILE2018} or the use of the supremizer stabilization method \cite{Ballarin2014,Rozza2007}, which consists into the enrichment of the velocity space by the usage of supremizer modes. These modes are computed such that a reduced version of the inf-sup condition is fulfilled. The latter approach usually employed in a finite element context has been also extended to a FV formulations \cite{STABILE2018}. In this work we rely on the supremizer stabilization method. After a proper enrichment of the POD velocity space it is possible to project the continuity equation onto the space spanned by the pressure modes giving rise to the following system:
\begin{eqnarray}\label{eq:reduced_system}
\left\{
\begin{matrix}
\nu \bm{B} \bm{a} - \bm{a}^T \bm{C} \bm{a}-\bm{H}\bm{b} =\bm{0},
\\
\bm{P} \bm{a}= \bm{0},
\end{matrix}
\right.
\end{eqnarray}
where the new matrix $\bm{P}$, is computed as follows:
\begin{align}
P_{ij}=\left(\chi_i, \nabla \cdot  \bm{\phi}_j \right)_{L^2(\Omega)}.
\end{align}
System \eqref{eq:reduced_system} can be solved respect to $\bm{a}$ and $\bm{b}$ in order to obtain the reduced order solution for velocity and pressure respectively. 

\subsection{Treatment of boundary conditions}\label{sec:boundary}\hspace*{\fill} \\
The current problem involves non-homogeneous Dirichlet boundary conditions at the inlet $\Gamma_{In}$. The term $\bm{f}(\bm{x},\bm{\mu})$ in \eqref{eq:navstokes} becomes $\bm{f}(\bm{x},\bm{\mu}) = (\mu_x,\mu_y) $ where $\mu_x$ and $\mu_y$ are the components of the velocity at $\Gamma_{In}$ along the $x$ and $y$ directions respectively. Non-homogeneous Dirichlet boundary conditions have been treated making use of the so called lifting control function. In this method the POD procedure is applied on a modified set of snapshots which have been homogenized in the following way:
\begin{equation}\label{eq:lifting}
\bm{u_i'} = \bm{u_i} - \mu_x \bm{{\phi_c}_x} - \mu_y \bm{{\phi_c}_y}, \quad \text{for} \quad i=1,...,N_s,
\end{equation}
where $\bm{{\phi_c}_x}$ and $\bm{{\phi_c}_y}$ are the two lifting functions which have at the inlet $\Gamma_{In}$ the following values $(1,0)$ and $(0,1)$ respectively. The approach used in this work to obtain the lifting functions involves solving two linear potential flow problems with the initial boundary conditions at $\Gamma_{In}$ being $(1,0)$ and $(0,1)$ respectively for $\bm{{\phi_c}_x}$ and $\bm{{\phi_c}_y}$.\par
The POD is then applied to the snapshots matrix $\bm{\mathcal{U'}} = [\bm{u'_1},\bm{u'_2},\dots,\bm{u'_{N_s}}]$ that contains only snapshots with homogeneous boundary conditions. It has to be noted that the way the lifting functions have been computed assures that they are divergence free and thus the new set of snapshots has the same property.\par
At the reduced order level, it is then possible to deal with any boundary velocity at $\Gamma_{In}$ (obviously the results will be more accurate if the prescribed velocity values are sufficiently close to the ones used during the training stage). If the new sample $\mu^\star$ (which has the new boundary velocity) is introduced in the online stage, one can compute the reduced velocity field as follows:
\begin{equation}
\bm{u}(\bm{x},\bm{\mu^\star})\approx \mu^\star_x \bm{{\phi_c}_x} + \mu^\star_y \bm{{\phi_c}_y}+\sum_{i=1}^{N_u} a_i(\bm{\mu^\star}) \bm{\phi_i} (\bm{x}).
\end{equation}
\section{Non-Intrusive PCE}\label{sec:PCE}

{\color{black}According to Polynomial Chaos (PC) theory which was formulated by Wiener \cite{wiener1938homogeneous}, real-valued multivariate Random Variables (RVs), such as the one considered in this work (the lift coefficient $C_l$) can be decomposed into an infinite sum of separable deterministic coefficients and orthogonal polynomials \cite{janya2017framework}. These polynomials are stochastic terms which depend on some mutually orthogonal Gaussian random variables. Once applied to our output of interest --- the lift coefficient $C_l$ --- such decomposition assumption reads  
\begin{equation}\label{eq:e1}
{C_l}^\star (\bm{\zeta}) = \sum_{i=0}^{\infty} {C_l}_i \psi_i (\bm{\zeta}).
\end{equation}
Here the random variable $\psi = (\alpha,U)$ is used to express the uncertainty in the angle of attack and inflow velocity. $ \psi_i (\bm{\zeta})$ is the $i^{th}$ polynomial and ${C_l}_i $ is the so-called $i^{th}$ stochastic mode. In practical application this series is truncated and only its first $P+1$ values are computed, namely
\begin{equation}\label{eq:e2}
{C_l}^\star (\bm{\zeta}) = \sum_{i=0}^{P} {C_l}_i \psi_i (\bm{\zeta}).
\end{equation}
In this work the orthogonal polynomial are called Hermite polynomials. These polynomials form an orthogonal set of basis functions in terms of Gaussian distribution \cite{ghanem2003stochastic}. In \eqref{eq:e2} $P+1$ is the number of Hermite polynomials used in the expansion and has to depend on the order of the polynomials chosen and the on dimension $n$ of the random variable vector $\bm{\zeta} = \{\zeta_1,\dots,\zeta_n\}$. More specifically, in an $n$-dimensional space, the number $P$ of
Hermite polynomials of degree $p$ is given by $ P+1 = \frac{(p+n)!}{p! n!}$\cite{ghanem2003stochastic}.}\par

\subsection*{Coefficients Computation}
The {\color{black} estimation} of the coefficients $ {C_l}_i (x)$ in \eqref{eq:e1} can be carried out in different ways. Among others, we mention the sampling based method and the quadrature method. The one here used is based on the sampling approach, following the methodology proposed by \cite{hosder2006non}. The coefficient calculation algorithm starts from a discretized version of equation \eqref{eq:e2}, namely
\[
\begin{bmatrix} 
    {C_l}_0^* \\
    {C_l}_1^* \\
    \vdots \\
    {C_l}_N^* 
    \end{bmatrix} = \begin{bmatrix} 
    \psi_1({\bm{\zeta}}_0) & \psi_2({\bm{\zeta}}_0) & \dots & \psi_P({\bm{\zeta}}_0) \\
    \psi_1({\bm{\zeta}}_1) & \psi_2({\bm{\zeta}}_1) & \dots & \psi_P({\bm{\zeta}}_1) \\
    \vdots & \ddots & \\
    \psi_1({\bm{\zeta}}_N) & \psi_2({\bm{\zeta}}_N) & \dots & \psi_P({\bm{\zeta}}_N) 
    \end{bmatrix} \begin{bmatrix} 
    {C_l}_0 \\
    {C_l}_1 \\
    \vdots \\
    {C_l}_P
    \end{bmatrix},
\] \label{eq:PCESystem}
where $N$ is the number of the samples taken. If $N$ coincides with the number of Hermite polynomials $P+1$ needed for the PCE expansion, the system above presents a square matrix and can be solved to determine the coefficients ${C_l}_i$ from the known output coefficients ${C_l}^\star_i$. In the most common practice, a redundant number of samples are considered and the system is solved in a least squares sense, namely
\begin{equation}
\bm{{C_l}} = (\bm{{L}}^T\bm{{L}})^{-1} \bm{{L}}^T \bm{{C_l}}^*,
\end{equation}
where $\bm{{L}}$, $ \bm{{C_l}}$ and $ \bm{{C_l}}^*$ denote the rectangular matrix in \eqref{eq:PCESystem}, the PCE coefficients vector and output vector, respectively.  
\vspace{-0.122cm}

\section{Numerical Results}\label{sec:results}
This section presents the results for the simulations carried out with the POD-Galerkin ROM and  PCE for UQ described in the previous sections. The first part of the analysis will be focused on the results obtained with the POD-Galerkin ROM. In the second part we will assess the performance of the  UQ technique on the airfoil problem, both when FOM and ROM simulation results are used to feed the PCE algorithm. The overall objective of the present section is in fact twofold. The first aim is to understand the influence of the samples distribution used to train the ROM in the results of the POD-Galerkin ROM. The second aim is to compare between the PCE UQ results obtained using full order model to those obtained with POD Galerkin-ROM. \par


\subsection{ROM results}\label{sec:ROM_res}
\begin{figure}
{
  \centering
  \begin{minipage}[b]{0.36\linewidth}
    \centering
    \includegraphics[width=\linewidth]{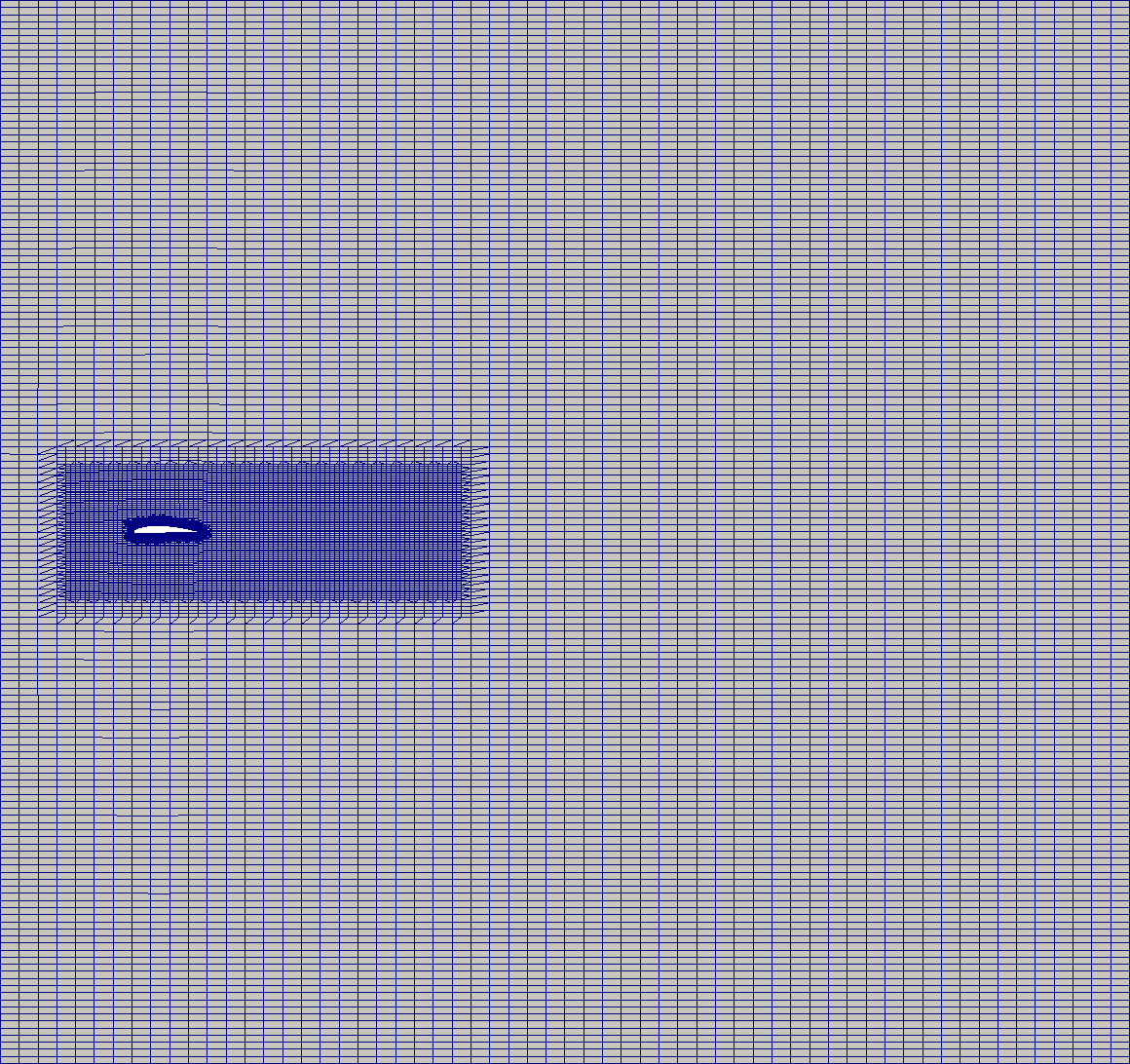}\label{fig:mesh} 
    \scriptsize(a)
    \end{minipage} 
  \begin{minipage}[b]{0.62\linewidth}
    \centering
    \includegraphics[width=\linewidth]{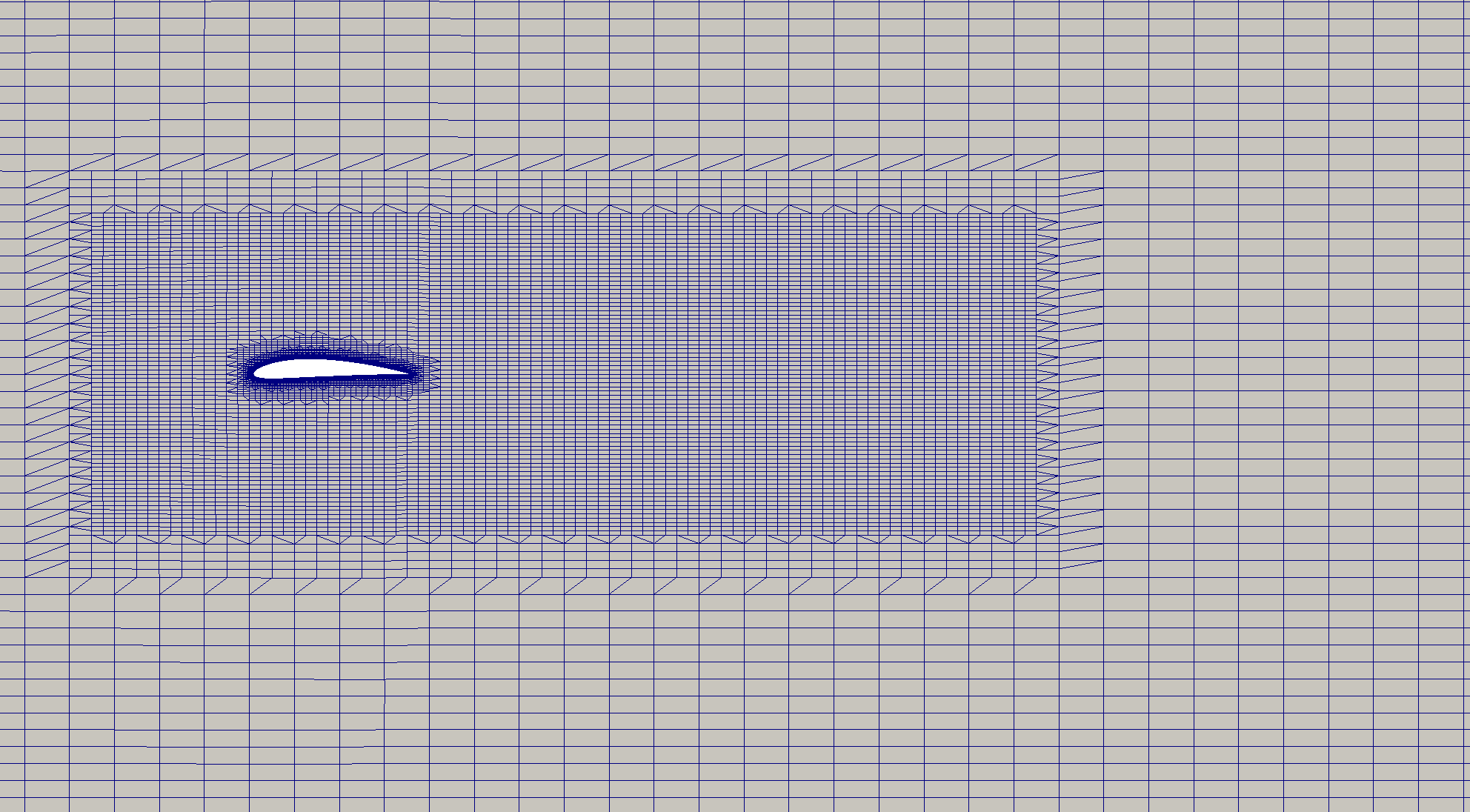}\label{fig:meshzoom}
        \scriptsize(b)
        \end{minipage} 
  \vspace{-0.2cm}\caption{(a) The OpenFOAM mesh used in the simulations. (b) A picture of the mesh zoomed near the airfoil.}\label{fig:Mesh} 
}
\end{figure}

\hspace*{\fill}\\
The FOM model used to generate the POD snapshots has been set up as reported in section \ref{sec:HF_model}. Making use of the computational grid shown in \autoref{fig:Mesh}, a set of simulations was carried out, selecting a Gauss linear numerical scheme for the approximation of gradients and Laplacian terms, and a bounded Gauss upwind scheme for the convective term approximation.\par

\begin{figure}
\centering
\includegraphics[width=0.8\textwidth]{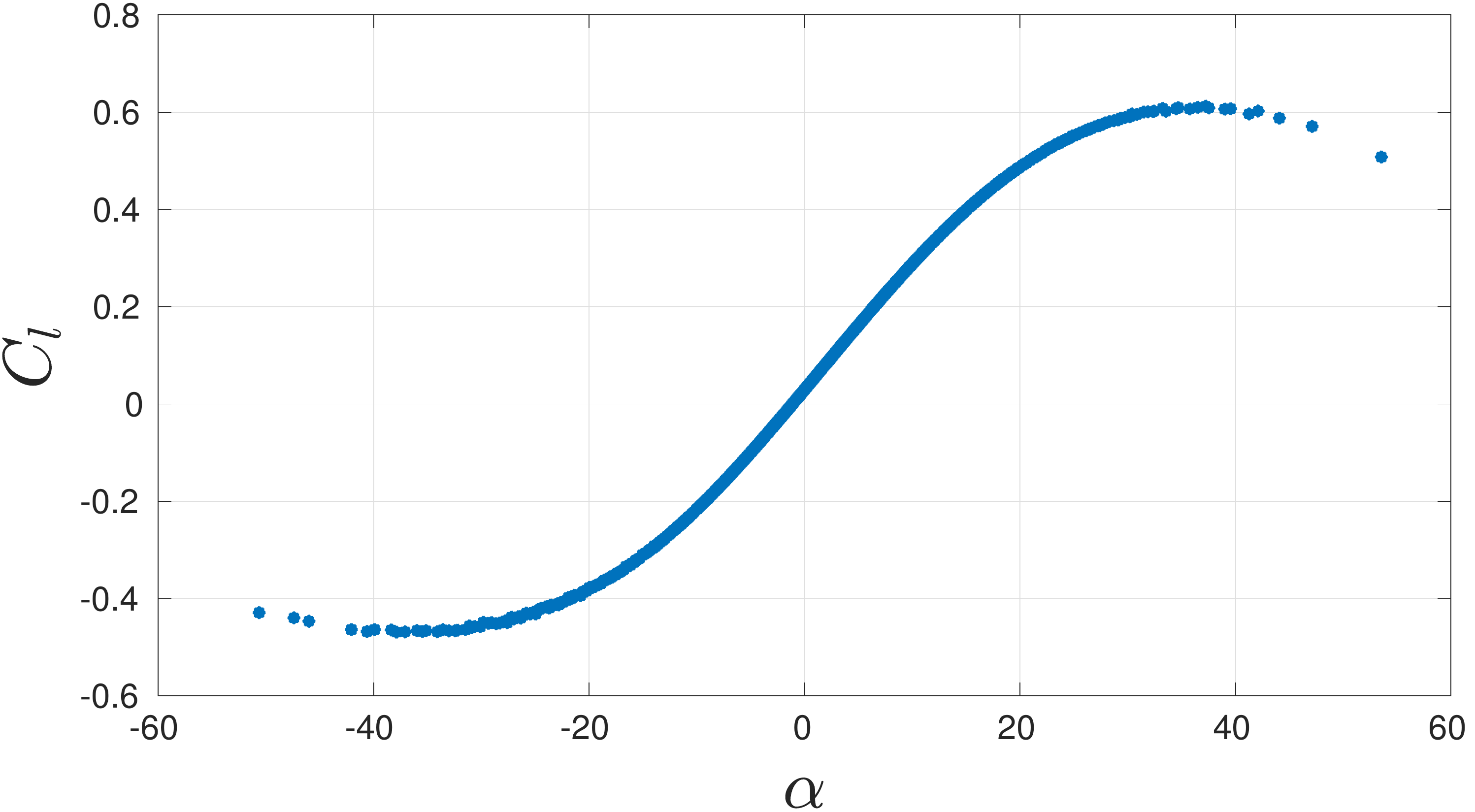}
\caption{The FOM lift coefficient for the first case, as it can be seen in the plot the vast majority of the samples is clustered around $\alpha = 0$.}\label{fig:LiftHFCase1}
\end{figure}

As mentioned, the parameters considered in the ROM investigation are the airfoil angle of attack and the magnitude of the inflow velocity at the inlet. The training of the POD-Galerkin ROM requires a suitable amount of snapshots (FOM solutions) to be available. Thus, $520$ samples have been produced and a single FOM simulation is launched for each sampling point. As for the distribution in the parameters space, the samples are obtained making use of the Latin Hyper Cube (LHC) \cite{Stein1987} sampling algorithm. \autoref{fig:LiftHFCase1} depicts the lift coefficient against the angle of attack curve obtained from a first FOM simulation campaign in which the $520$ samples were generated imposing mean values of  $100\,\text{m}/\text{s}$ and $0^{\circ}$ and variances of $20\,\text{m}/\text{s}$ and $300$, for velocity and angle of attack, respectively. {\color{black}As can be appreciated in the picture, the lift coefficient values do not significantly depend on the inflow velocity. In fact sampling points with equal $\alpha$ and different $U_{\infty}$ values, result in practically identical output. For this reason, the input-output relationship appears like a curve in the $C_l$-$\alpha$ plane. We also point out that this is a consequence of considering a nondimensionalized force a $C_l$ as our output, rather then the corresponding dimensional lift values.}\par 

\begin{figure}
\centering
\includegraphics[width=0.8\linewidth]{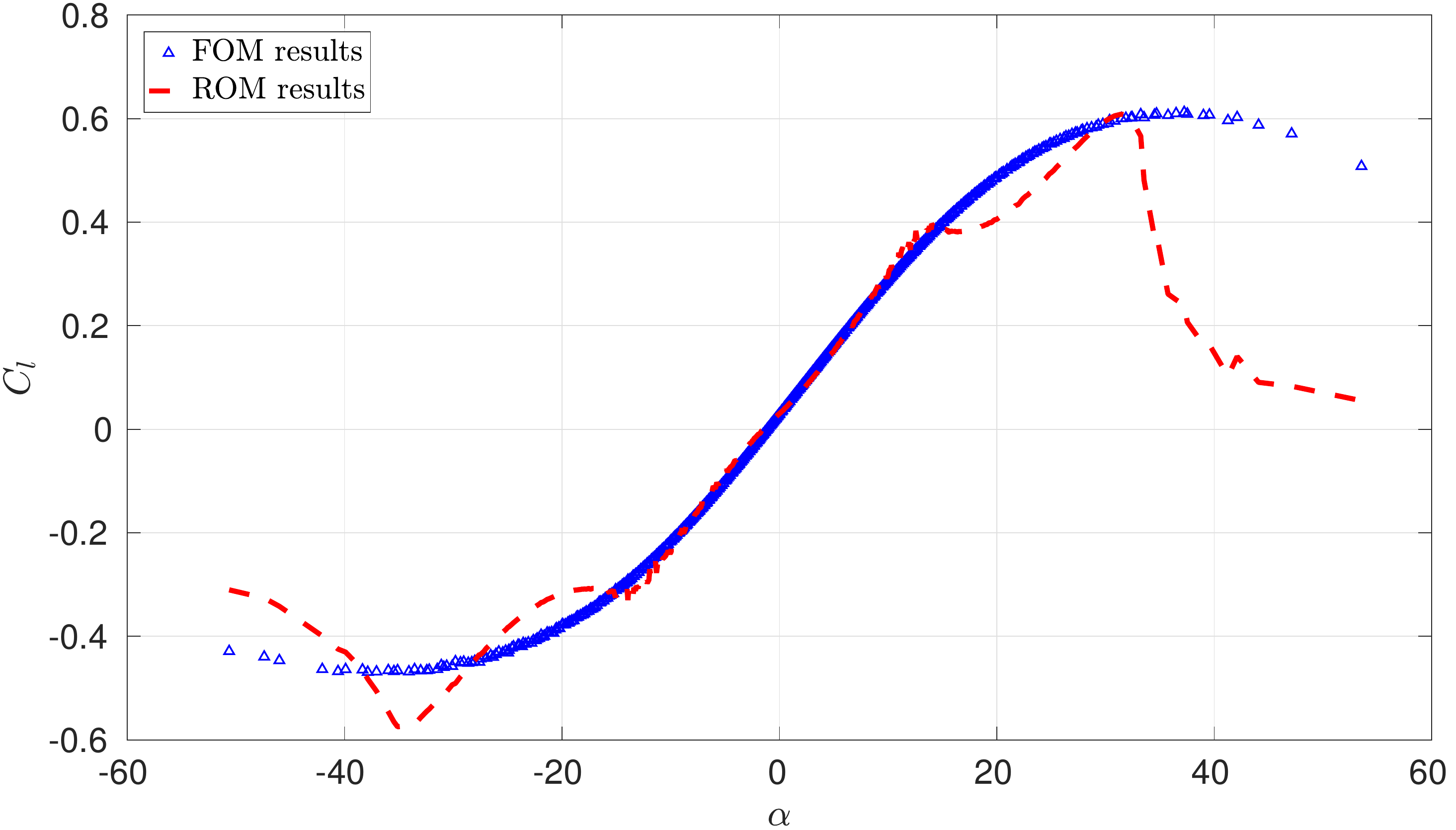} 
\caption{The first sampling case : the full order lift coefficients curve versus the ROM reconstructed one with $10$ modes used for each of velocity, pressure and supremizer fields.}    \label{fig:LiftPOD_case1a}  
\end{figure}

\begin{figure}
\centering
\includegraphics[width=0.8\linewidth]{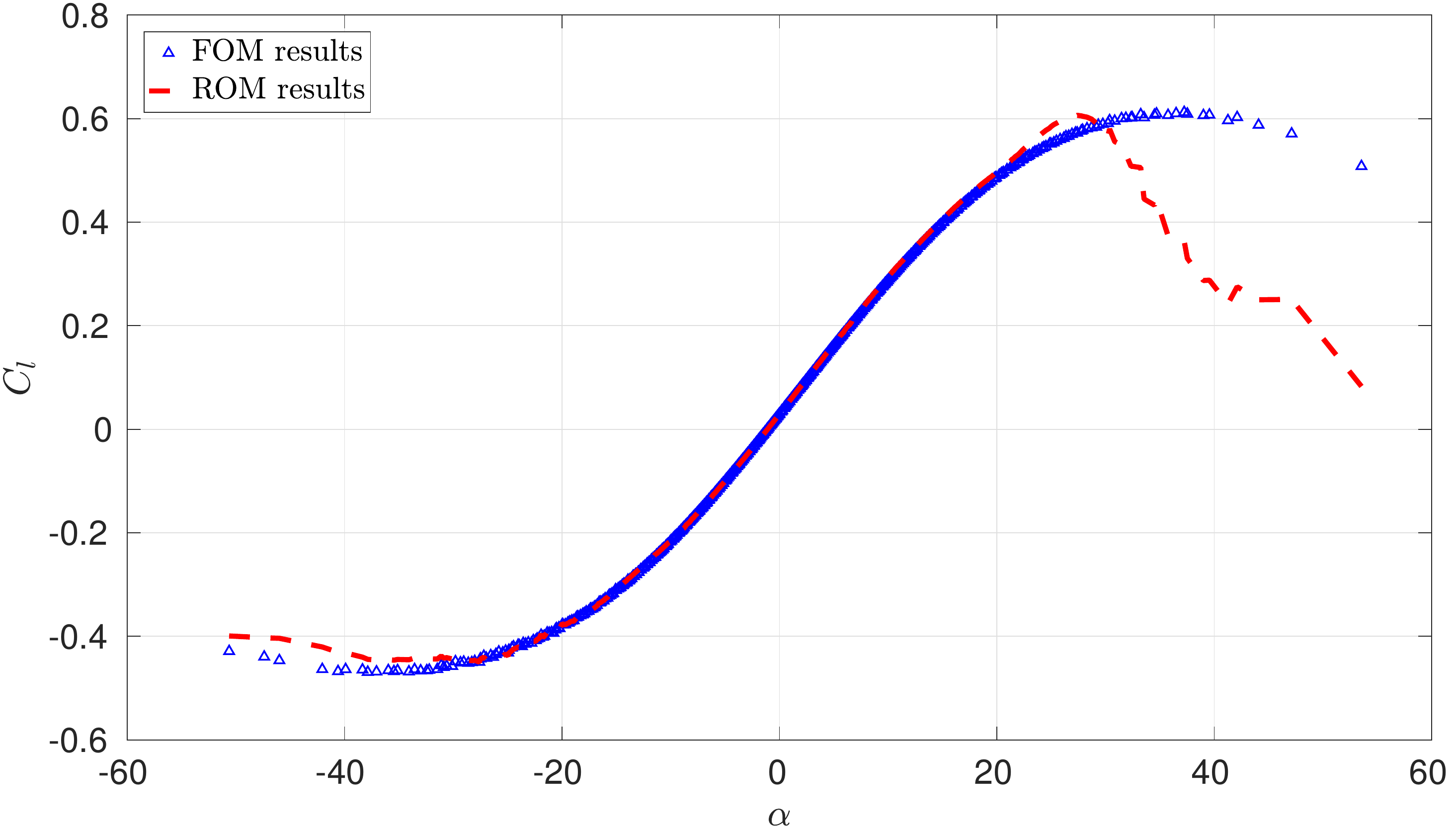}
\caption{The first sampling case : the full order lift coefficients curve versus the ROM reconstructed one with $30$, $10$ and $10$ modes are used for velocity, pressure and supremizer fields respectively.}\label{fig:LiftPOD_case1b} 
\end{figure}

The POD modes are generated after applying POD onto the snapshots matrices of the flow fields obtained in the simulation campaign. After such offline phase, the computation of the reduced order fields is performed in the online stage, as presented in \autoref{sec:ROM}. In this first reproduction test, we performed a single reduced simulation in correspondence with the velocity and the angle of attack used to generate each offline snapshot. This means that we used the same sample values both in the online stage and in the offline one. The ROM results of the reproduction test for the lift coefficient are reported in \autoref{fig:LiftPOD_case1a}. The figure refers to the ROM results obtained considering $10$ modes for the discretization of velocity, pressure and supremizer fields. The plot shows that the reconstruction of the lift coefficient is only accurate in the central region. In the lateral regions the lift coefficient computed with the ROM solution does not match the corresponding FOM solution. The poor quality of the ROM prediction on the lift coefficient, as well as of the forces acting on the airfoil, is a direct consequence of the fact that the fields were not reconstructed in an accurate way. For the particular physical phenomenon this inaccuracy may be even more undesirable, since the stall occurs in these regions. One might originally guess that the problem can be mitigated by increasing the amount of POD modes. Yet, as \autoref{fig:LiftPOD_case1b} clearly shows, even increasing the modes for velocity to $30$ is not solving the problem. In this particular case adding more modes will not solve the problem since the energy added by considering more modes is negligible. An explanation of the poor performances of the ROM model in the stall region may be instead associated to the distribution of the offline samples used to generate the POD snapshots. In fact, the samples generated with LHC are distributed around a mean value of the angle of attack of $0^{\circ}$ and their density is rather coarse in the stall regions.\par
\begin{figure}
\centering
\includegraphics[width=0.8\textwidth]{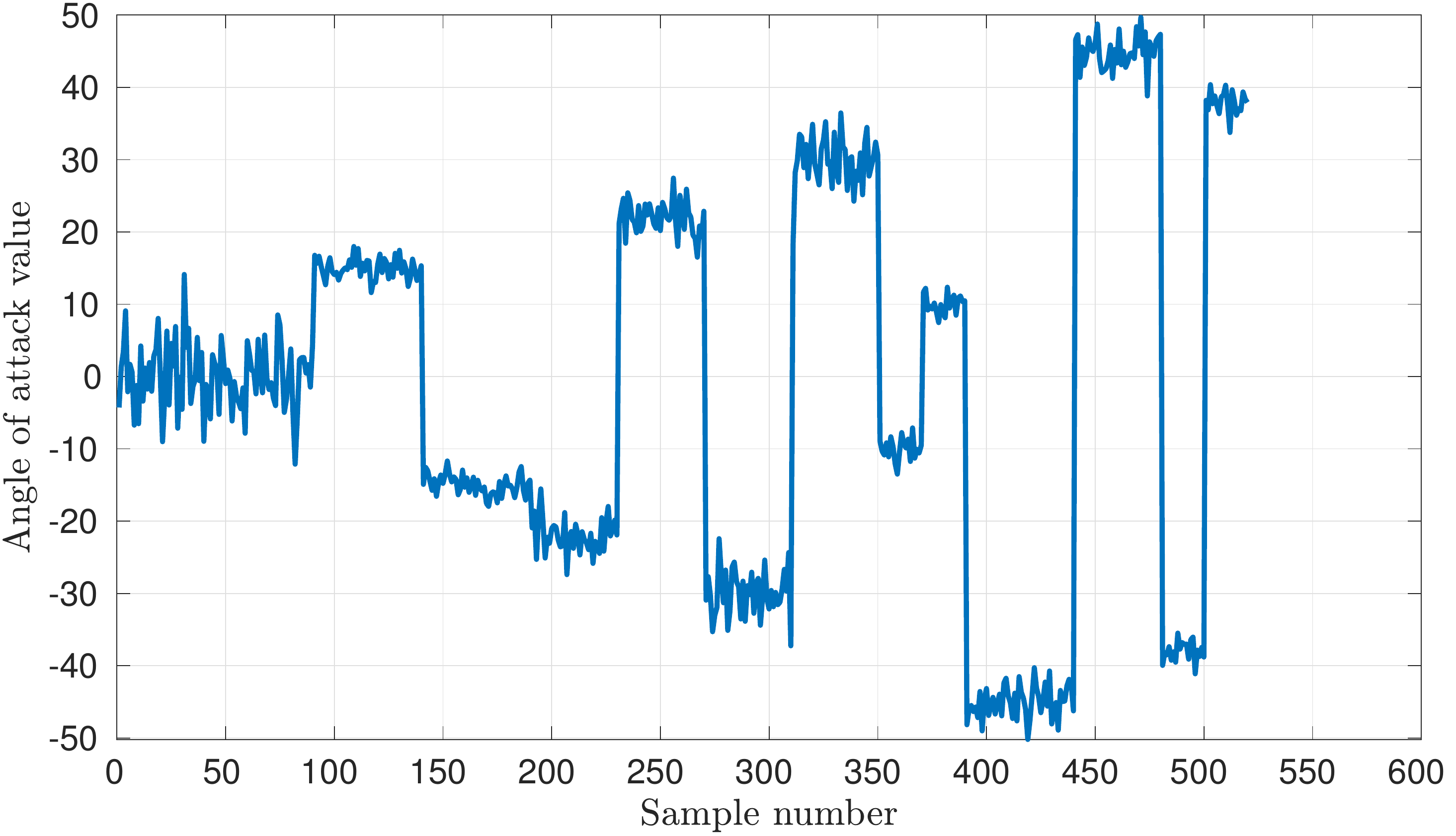}
\caption{The angle of attack samples for the second case.}\label{fig:sampang1}
\end{figure}

\begin{figure}
\centering
\includegraphics[width=0.8\textwidth]{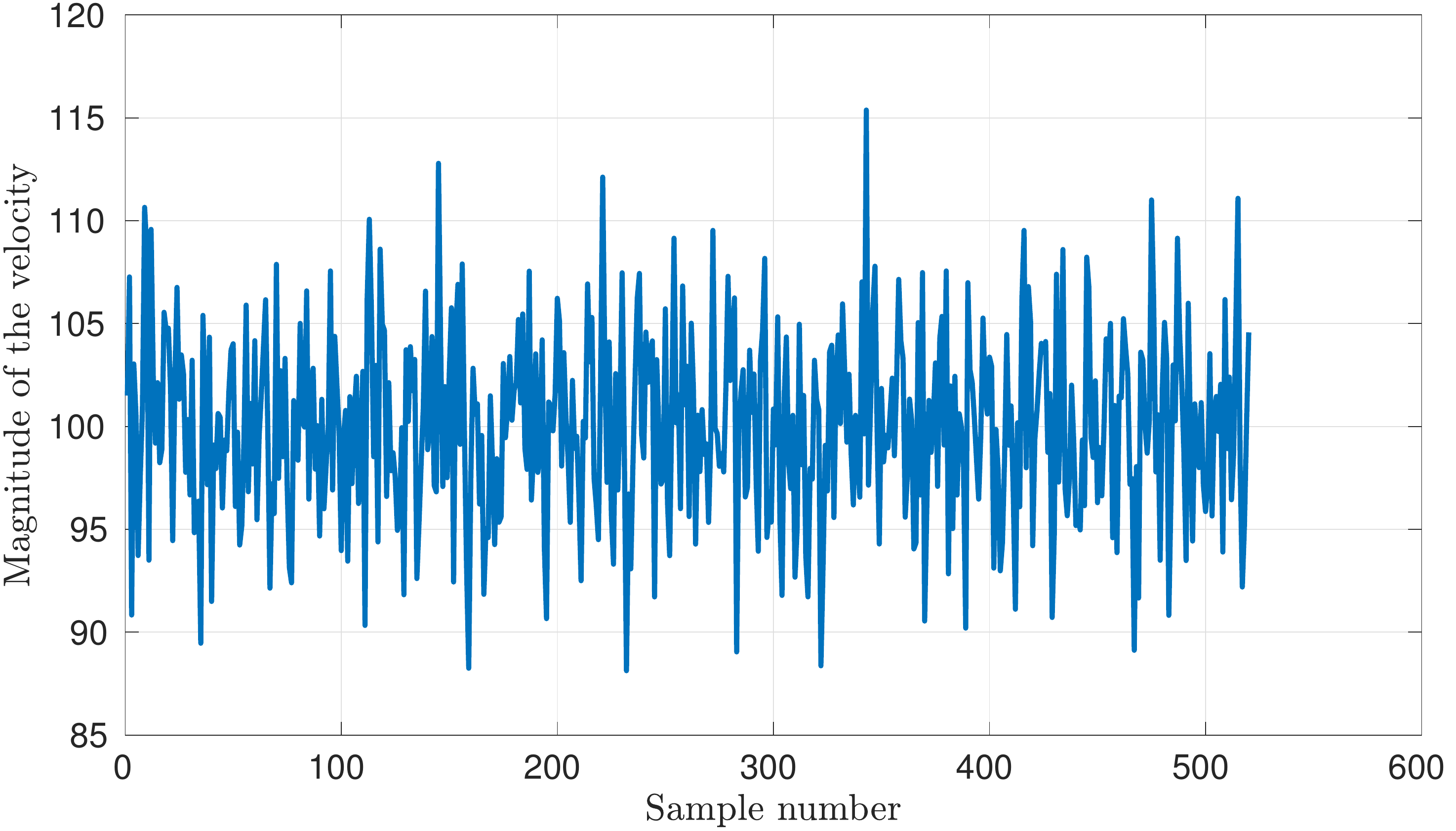}
\caption{The magnitude of velocity samples for the second case.}\label{fig:sampvel1}
\end{figure}

\begin{figure}
\centering
\includegraphics[width=0.8\textwidth]{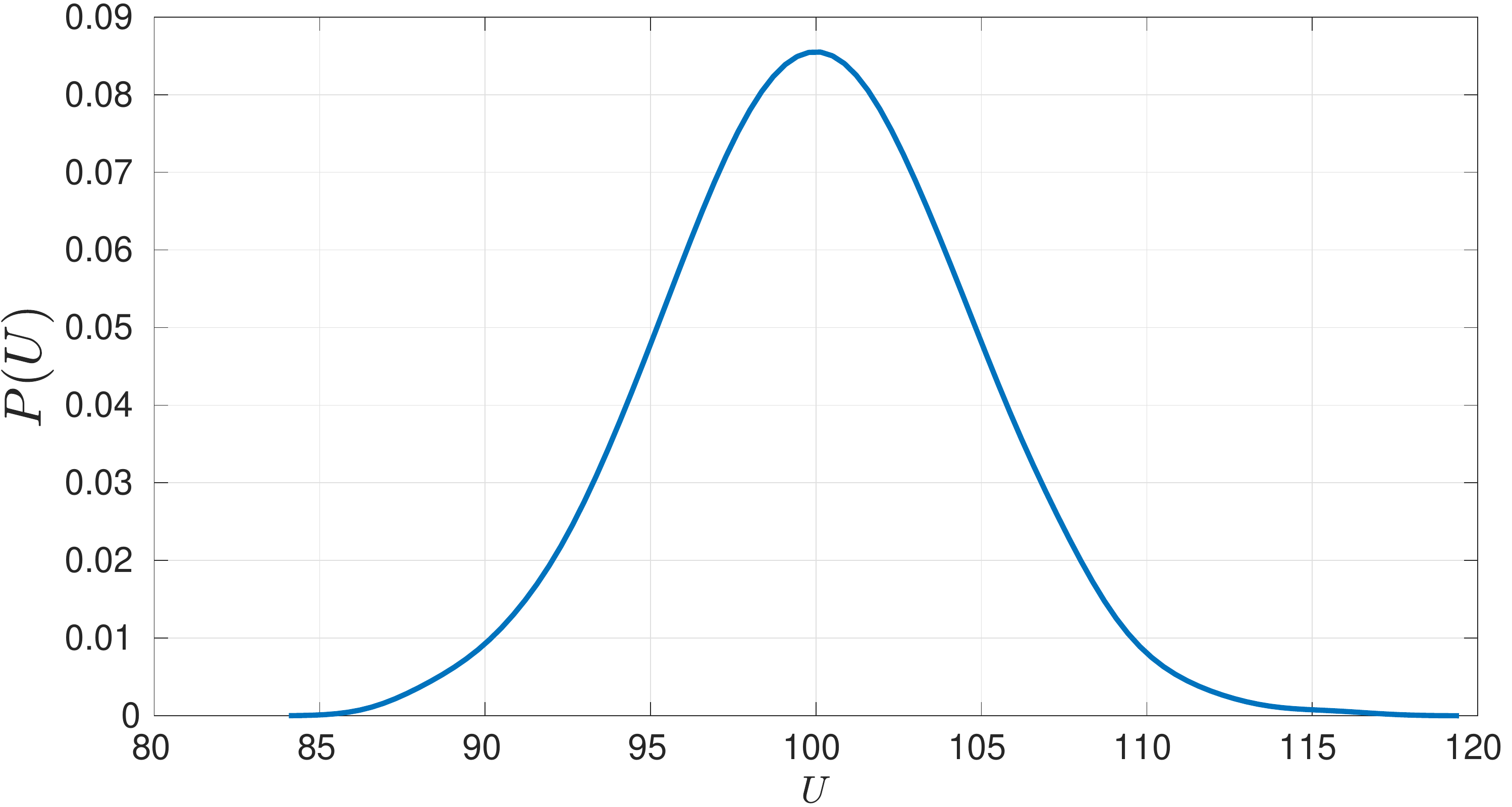}
\caption{The PDF of the input data set for the magnitude of the velocity.}\label{fig:pdf_input_vel}
\end{figure}

\begin{figure}
\centering
\includegraphics[width=0.8\textwidth]{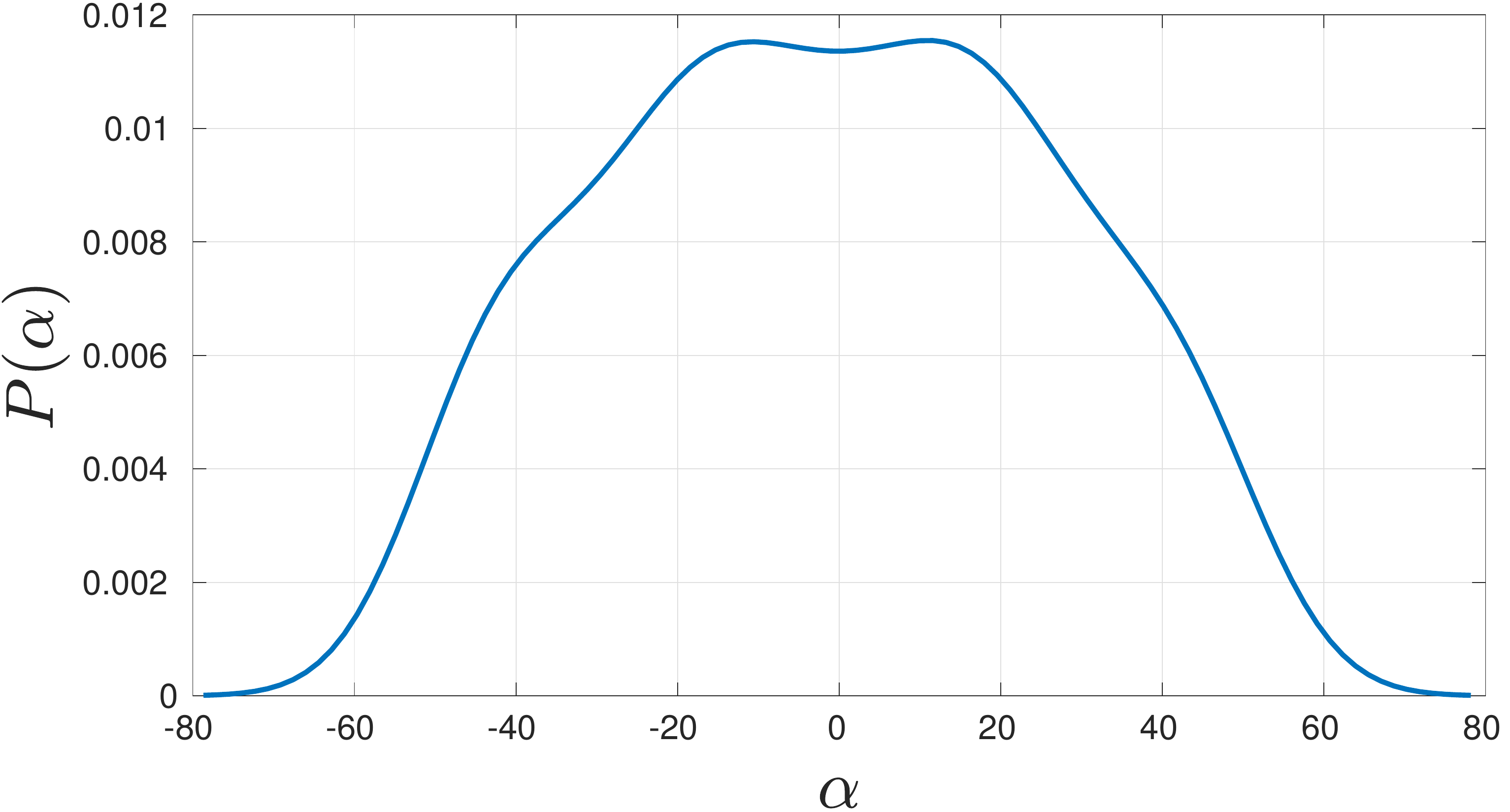}
\caption{The PDF of the input data set for the angle of attack.}\label{fig:pdf_input_ang}
\end{figure}

 
To confirm such deduction, we tried to generate the snapshots by means of a different set of samples generated so as to be more dense in the stall regions. More specifically, we have generated thirteen different groups of samples in which the velocity mean and variance were kept fixed at values of $100\,\text{m}/\text{s}$ and $20\,\text{m}/\text{s}$ respectively, while the mean and the variance for the angle of attack were varied in each set as summarized in \autoref{table:case2}. As illustrated in \autoref{fig:sampang1} and \autoref{fig:sampvel1} for angle of attack and velocity respectively, the overall training data for the ROM offline phase has been generated combining into a single sample set all the 13 groups generated by means of LHC algorithm. Finally, \autoref{fig:pdf_input_vel} and \autoref{fig:pdf_input_ang} depict the Probability Density Function (PDF) of the overall input parameters set.\par

\begin{table}[htp]
\centering
{
\begin{tabular}{ c | c | c | c  }
Group Number & N & $E(\alpha) \quad \text{in} \quad ^{\circ}$ & $\sigma(\alpha)  \quad \text{in} \quad ^{\circ}$ \\
\hline  			
$1$  & $90$  & $0$ & $20$ \\
$2$  & $20$  & $-10$ & $2$  \\
$3$  & $20$  & $10$ & $2$  \\
$4$  & $50$  & $-15$ & $2$  \\
$5$  & $50$  & $15$ & $2$  \\
$6$  & $40$  & $-22$ & $5$  \\
$7$  & $40$  & $22$ & $5$  \\
$8$  & $40$  & $-30$ & $10$ \\
$9$  & $40$  & $30$ & $10$  \\
$10$  & $20$  & $-38$ & $2$  \\
$11$  & $20$  & $38$ & $2$  \\
$12$  & $50$  & $-45$ & $5$  \\
$13$  & $40$  & $45$ & $5$  
\end{tabular}}\caption{{The mean and variance for the group of samples that form the training set in the second case.}}
\label{table:case2}
\end{table}
After running the offline phase, we applied POD-Galerkin reduced order approach on the new set of snapshots generated. \autoref{tab:cum_eig} shows the cumulative eigenvalues for the correlation matrices built by the snapshots obtained for velocity, pressure and supremizer fields. Only the values of the cumulative eigenvalues up to the fifteenth mode are listed. Yet, such data indicate that by using $15$ modes in the online phase, we can recover $99.9\%$ of the energy embedded in the system.\par
\begin{figure}
\centering
\includegraphics[width=0.8\textwidth]{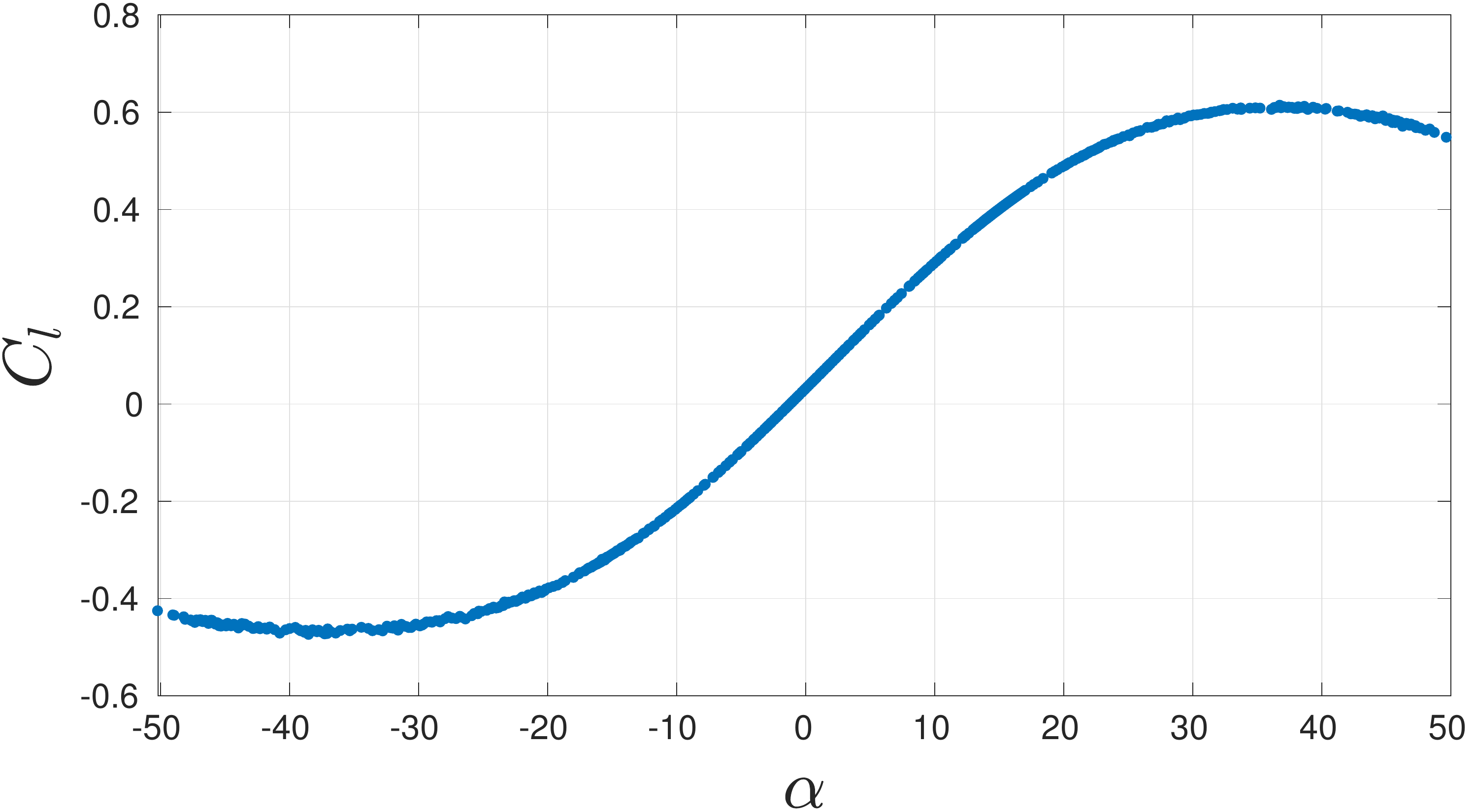}
\caption{The FOM lift coefficient as a function of the angle of attack $\alpha$ for the second case. Compared to \autoref{fig:LiftHFCase1}, the samples span a wider range with sufficient accuracy sought by the ROM reconstruction.}\label{fig:LiftHFCase2}
\end{figure}

\begin{figure}
{
  \centering
  \begin{minipage}[b]{0.8\linewidth}
    \centering
    \includegraphics[width=\linewidth]{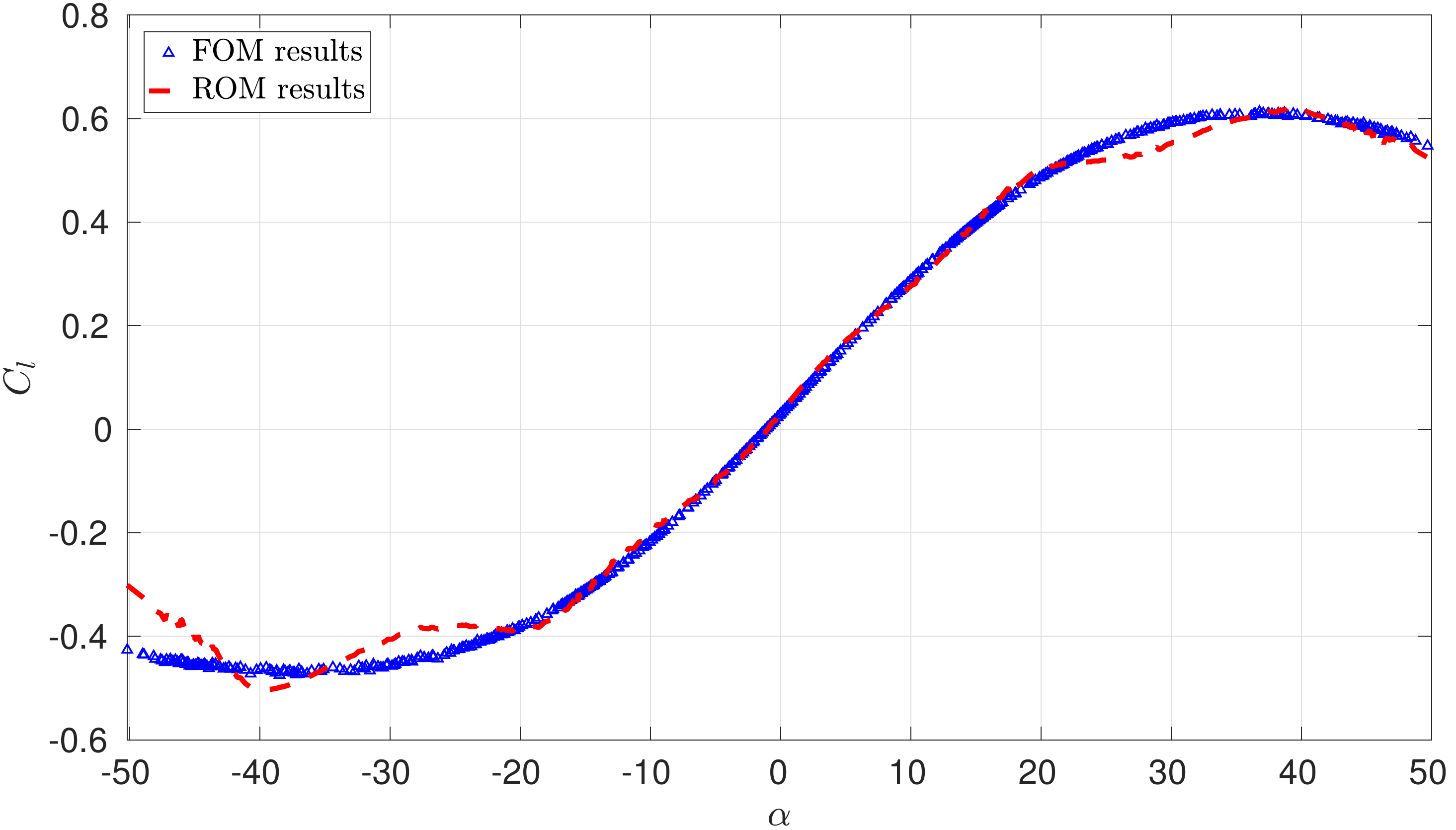} 
    \scriptsize(a)
 
    \end{minipage}
    \begin{minipage}[b]{0.8\linewidth}
    \centering
    \includegraphics[width=\linewidth]{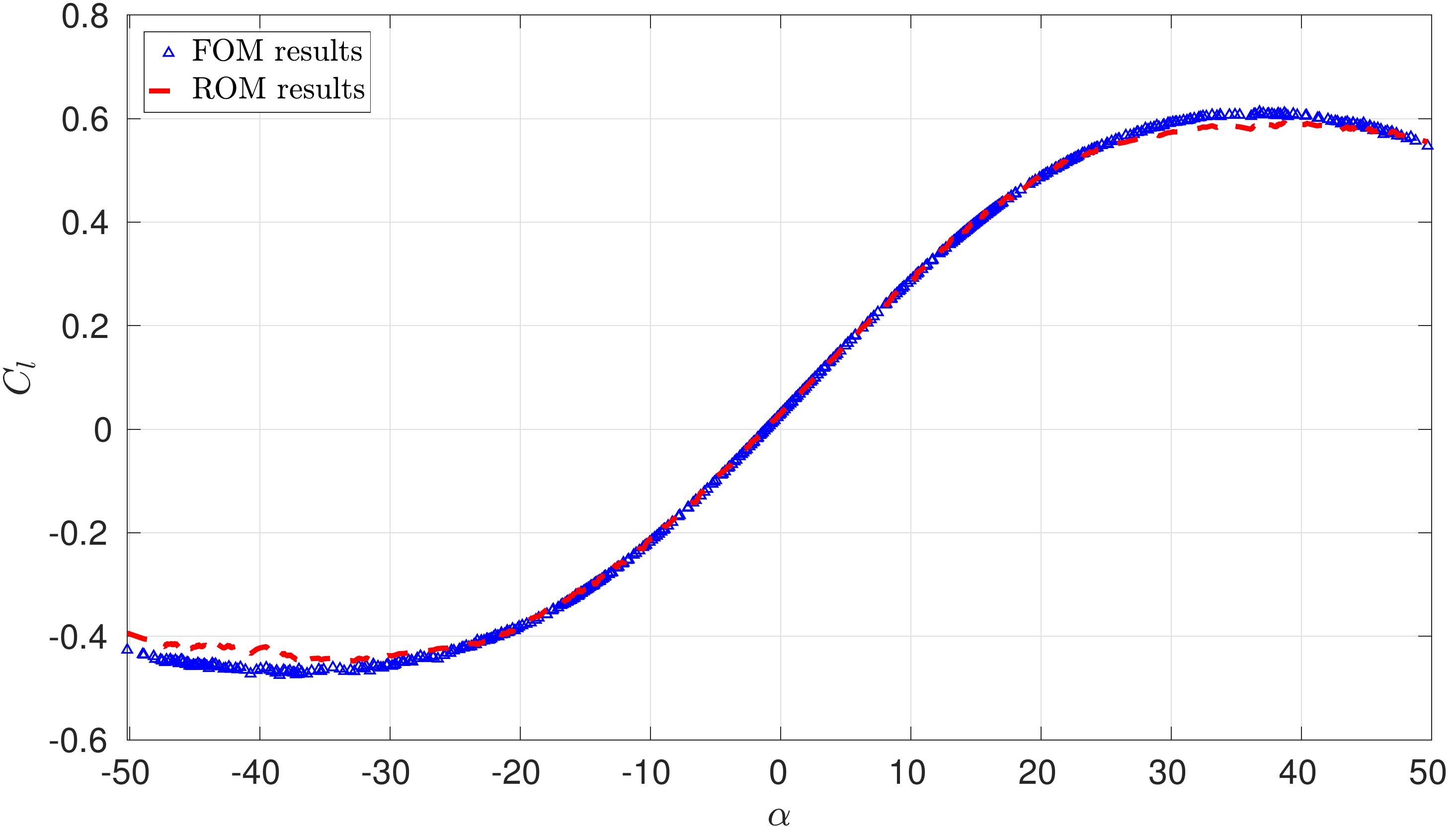}
        \scriptsize(b)
 
        \end{minipage} 
  \vspace{-0.2cm}\caption{A comparison between FOM and ROM reconstructed lift coefficients for the second case (a) $10$ modes are used for each of velocity, pressure and supremizer fields. (b) $20$, $10$ and $10$ modes are used for velocity, pressure and supremizer fields, respectively. In both figures we have online parameters set that coincide with the offline ones.}\label{fig:LiftROMcase2_1} 
}
\end{figure}

\begin{figure}
{
  \centering
  \begin{minipage}[b]{0.48\linewidth}
    \centering
    \includegraphics[width=\linewidth]{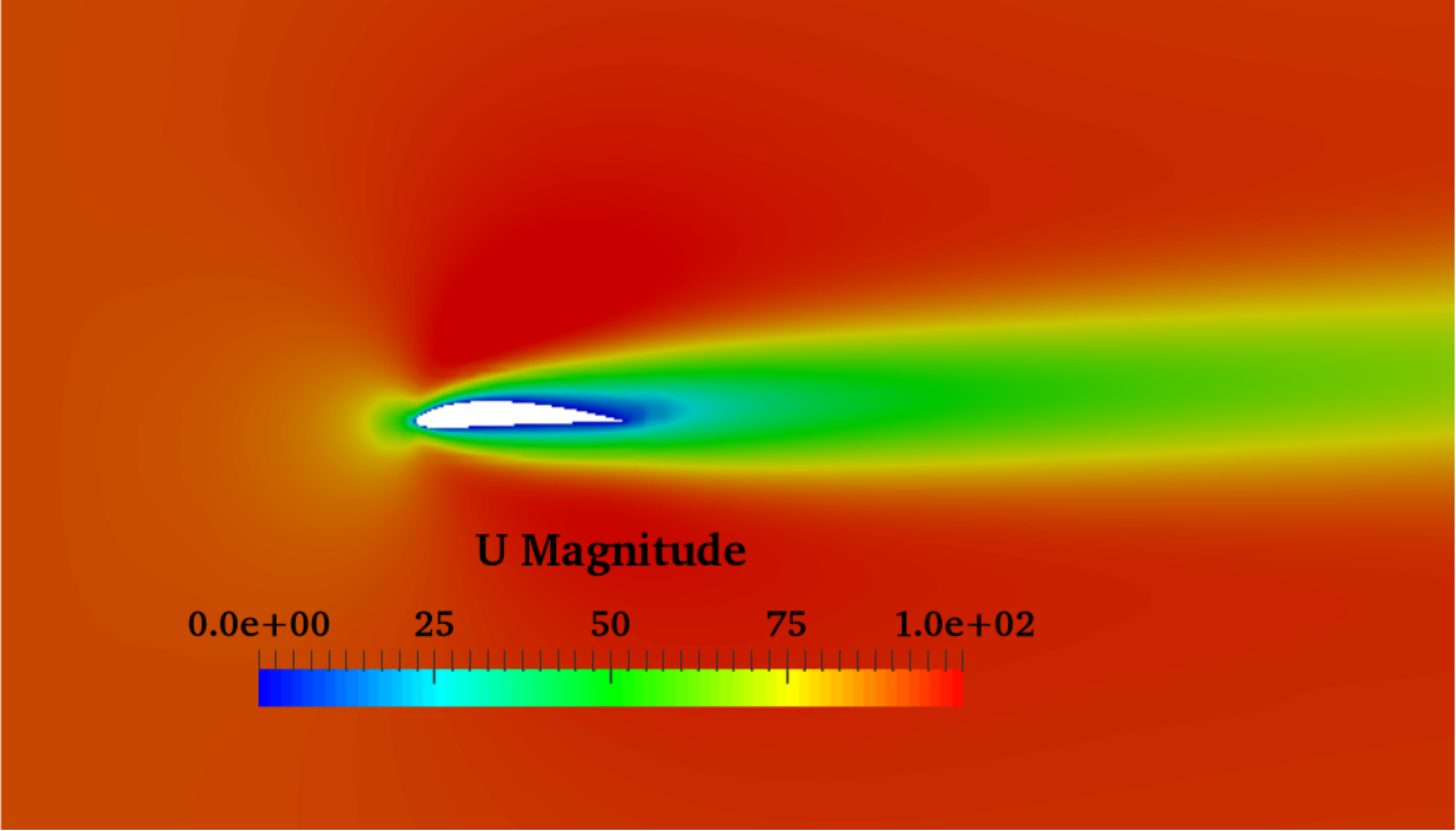} 
    \scriptsize(a)

    \end{minipage}
  \begin{minipage}[b]{0.48\linewidth}
    \centering
    \includegraphics[width=\linewidth]{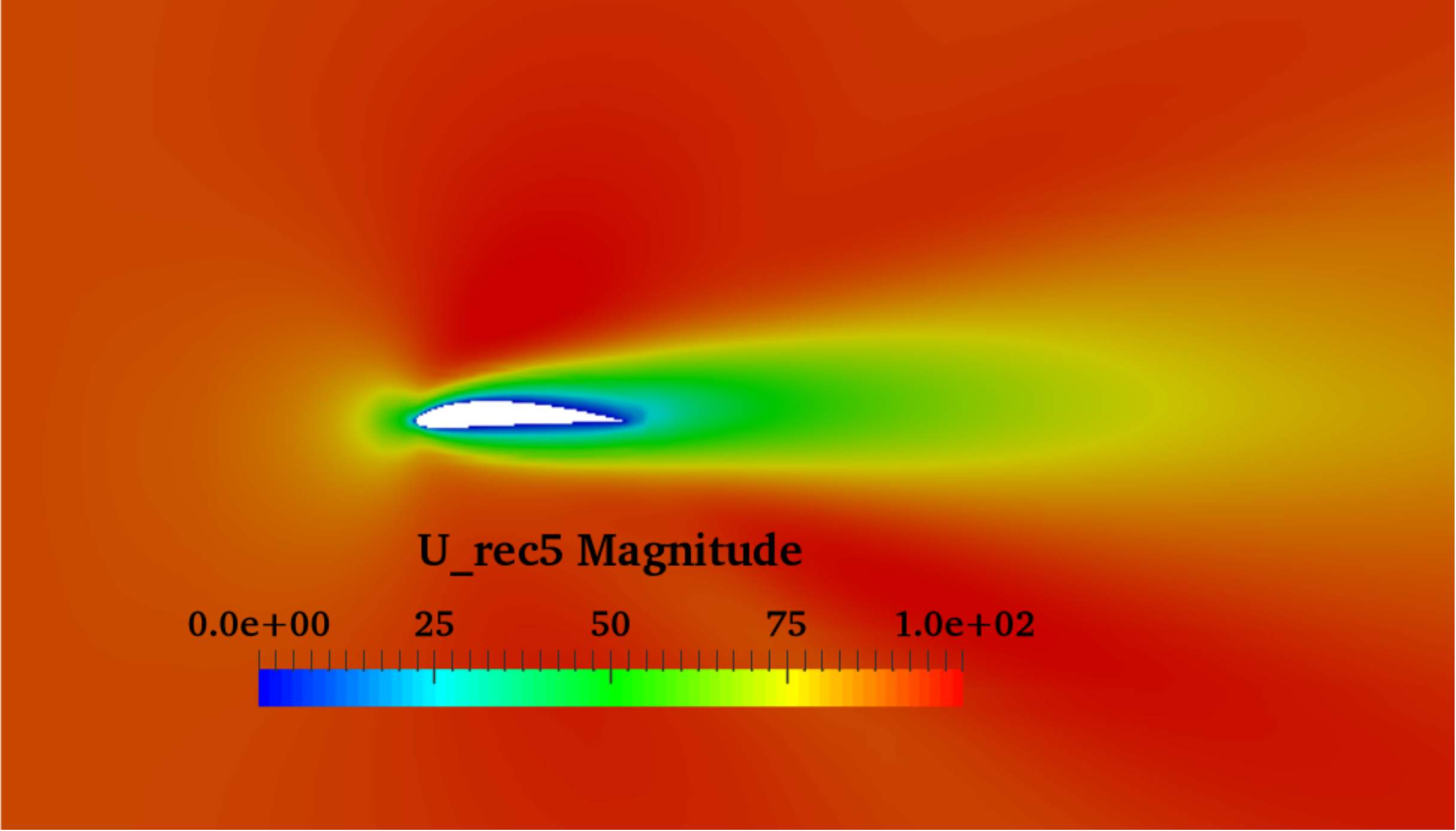}
        \scriptsize(b)
        \end{minipage} 
   \begin{minipage}[b]{0.48\linewidth}
    \centering
    \includegraphics[width=\linewidth]{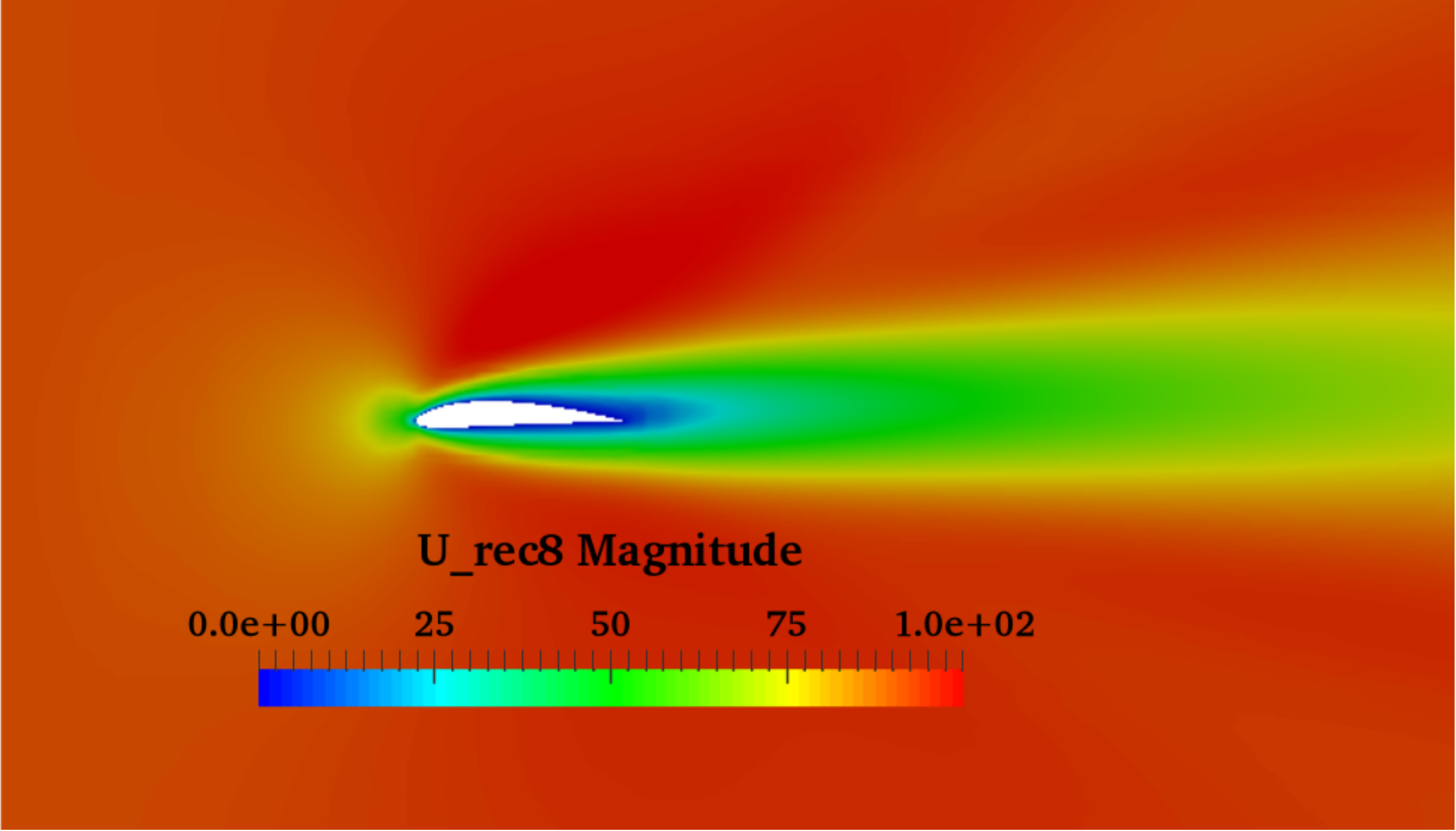}
        \scriptsize(c)
        \end{minipage} 
   \begin{minipage}[b]{0.48\linewidth}
    \centering
    \includegraphics[width=\linewidth]{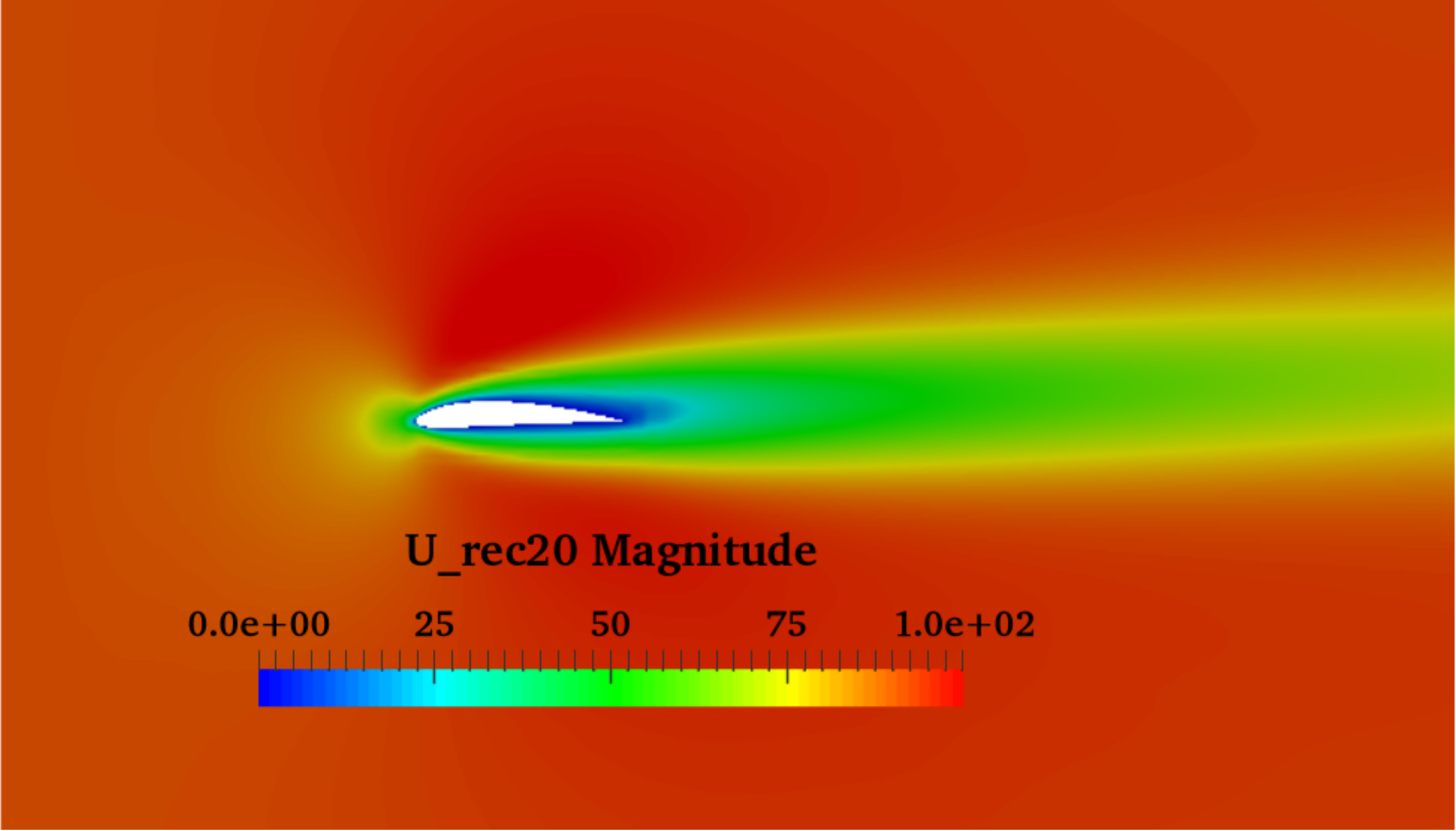}
        \scriptsize(d)
        \end{minipage} 
  \vspace{-0.2cm}\caption{The full order velocity field for the parameter $\bar{\bm{\mu}} = (90.669$ $m/s$ , $5.4439$ $m/s)$ and a comparison with the reconstructed field by means of different number of modes for velocity, pressure and supremizer fields. (a) FOM field (b) ROM velocity field with $3$, $10$ and $10$ modes used. (c) ROM velocity field with $8$, $10$ and $10$ modes used. (d) ROM velocity field with $20$, $10$ and $10$ modes used.}\label{fig:vel_comp} 
}
\end{figure}

 
\begin{table}[htp]
\centering
{
\begin{tabular}{ c | c | c | c  }
N Modes & $\bm{u}$ & $p$ & $\bm{u}_{sup}$\\
\hline  			
$1$ & $0.679635$ & $0.738828$ & $0.557189$ \\ 
$2$ & $0.930038$ & $0.960781$ & $0.987862$ \\ 
$3$ & $0.955239$ & $0.990746$ & $0.995984$ \\ 
$4$ & $0.971768$ & $0.998833$ & $0.999228$ \\ 
$5$ & $0.981370$ & $0.999730$ & $0.999796$ \\ 
$6$ & $0.987603$ & $0.999880$ & $0.999927$ \\ 
$7$ & $0.992311$ & $0.999945$ & $0.999975$ \\ 
$8$ & $0.994793$ & $0.999963$ & $0.999983$ \\ 
$9$ & $0.996651$ & $0.999976$ & $0.999990$ \\ 
$10$ & $0.997914$ & $0.999982$ & $0.999993$ \\ 
$11$ & $0.998679$ & $0.999987$ & $0.999995$ \\ 
$12$ & $0.999165$ & $0.999991$ & $0.999997$ \\ 
$13$ & $0.999492$ & $0.999993$ & $0.999998$ \\ 
$14$ & $0.999700$ & $0.999995$ & $0.999999$ \\ 
$15$ & $0.999806$ & $0.999996$ & $0.999999$
\end{tabular}}\caption{{Cumulative Eigenvalues of the correlation matrices for velocity, pressure and supremizer fields.}}
\label{tab:cum_eig}
\end{table}
\autoref{fig:LiftHFCase2} displays the FOM lift coefficients of the airfoil corresponding to the set of angles of attack previously introduced in \autoref{table:case2}. \autoref{fig:LiftROMcase2_1} shows the results obtained with the reduced order model trained with the samples summarized in \autoref{table:case2}. Here, the online phase was carried out using the same samples used in the offline stage. To provide a quantitative evaluation of the results, we used $L^2$ relative error, computed as follows:
\begin{equation}
\varepsilon = 100 \frac{\sqrt{\sum_{t=1}^n (C^{FOM}_{l_t} - C^{ROM}_{l_t})^2}}{\sqrt{\sum_{t=1}^n (C^{FOM}_{l_t})^2}}  \%,
\end{equation}
where $n$ is the number of sampling points, $C^{FOM}_{l_t}$ and $C^{ROM}_{l_t}$ are the $t-$th sample point lift coefficients for FOM and ROM, respectively. In this case we have a relative error of $6,69\%$ in $L^2$ norm between the FOM and the ROM lift coefficients, when $10$ modes have been used in the online phase for each of velocity, pressure and supremizer fields. Using 10 additional modes for velocity, results instead in a $3,75\%$ error. The corresponding plots in \autoref{fig:LiftROMcase2_1} clearly suggest that the qualitative behavior of the ROM lift output was substantially improved with respect to the first case. This improvement in the prediction of the ROM lift coefficients is due to a more accurate reproduction of the ROM fields. This is highlighted by \autoref{fig:vel_comp}, which shows the FOM velocity field along with different reconstructed surrogate fields obtained employing different number of modes at the projection stage.\par
%
%
\subsection{PCE results}\label{sec:PCE_res}
\begin{figure}
{
  \centering
  \begin{minipage}[b]{\textwidth}
    \centering
    \includegraphics[width=0.75\textwidth]{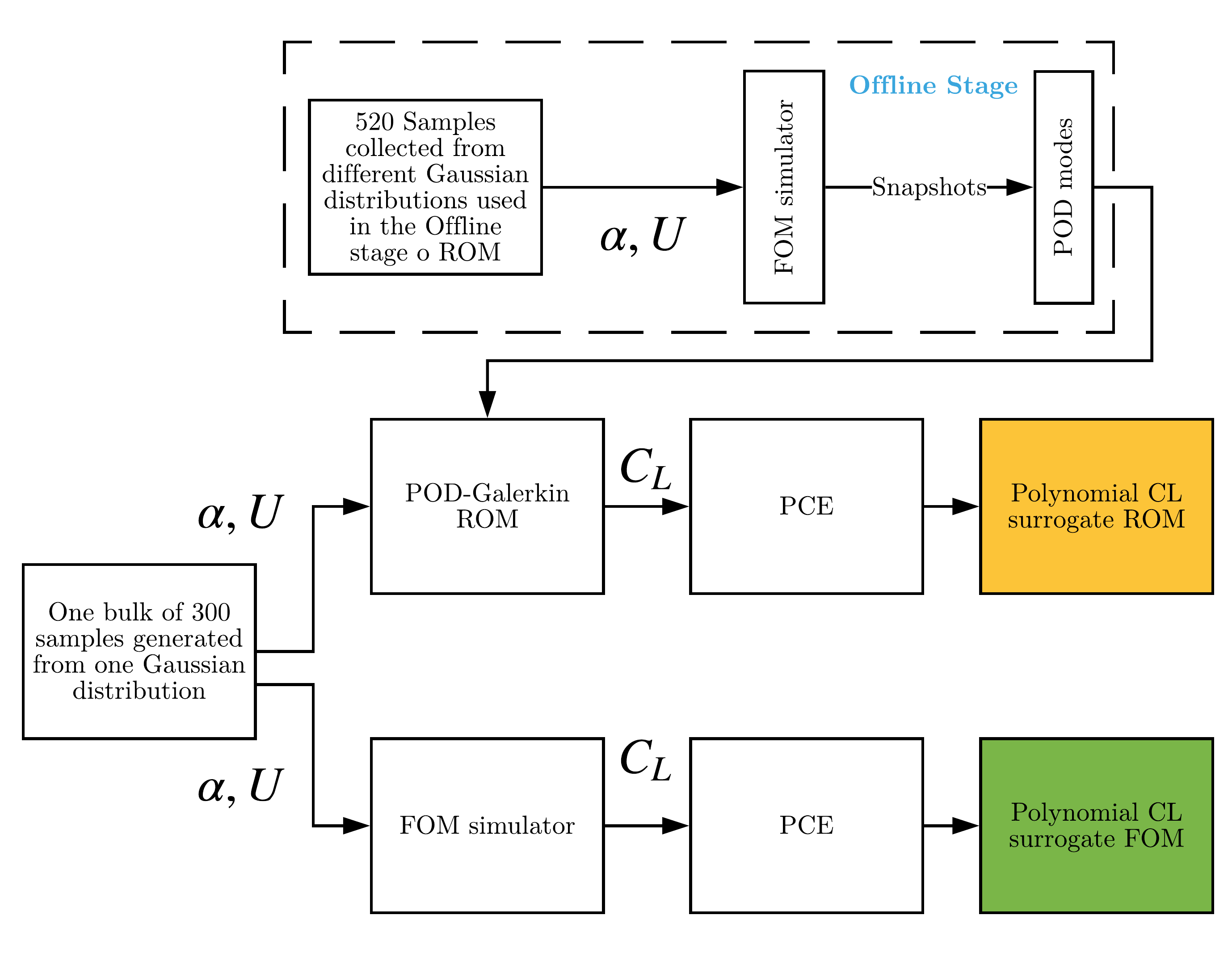} 
    \scriptsize(a)
    \end{minipage} 
  \begin{minipage}[b]{\textwidth}
    \centering
    \includegraphics[width=0.75\textwidth]{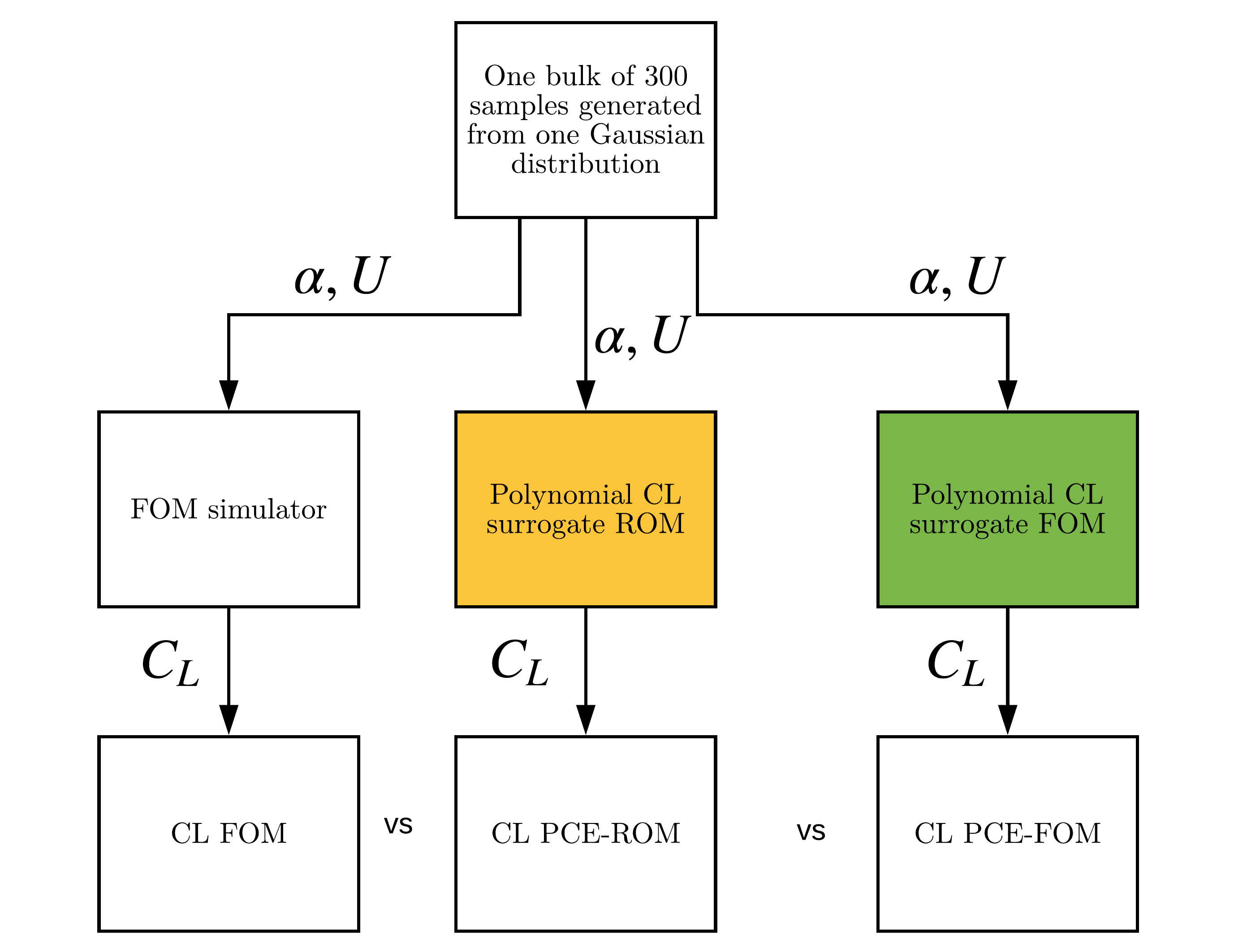}
        \scriptsize(b)
        \end{minipage} 
  \vspace{-0.2cm}\caption{The flowcharts describing procedure followed in the numerical simulations for the UQ model generation and validation campaign respectively.
                          The top scheme focuses on the procedure adopted for the generation of the UQ model, and in particular on the identification of the PCE coefficients. The polynomial surrogate based
                          on the full order model (indicated in green) has been generated using 300 Gaussian distributed samples in the $\alpha,U$ space. The same samples have been used to obtain the
                          polynomial surrogate input-output relationship for the POD-Galerkin ROM (denoted by the yellow box). Note that the ROM used in this simulation campaign has been trained by means
                          of 520 samples in the $\alpha,U$ space, organized in 13 Gaussian distributed bulks, as reported in Table \ref{table:case2}. Finally, the bottom flowchart illustrates
                          the PCE validation campaign. Here, 300 sample points in the input space have been used to obtain the corresponding output with the full order model, with the polynomial
                          UQ surrogate trained with the FOM simulations (green box), and with the polynomial UQ surrogate trained with the ROM simulations (yellow box)}\label{fig:pce_flowchart}
}
\end{figure}

\hspace*{\fill}\\
The aim of the present section is to evaluate the performance of the PCE algorithm implemented for the fluid dynamic problem at hand. To better describe the amount of simulations carried out to both train and validate the UQ PCE model implemented, we present in \autoref{fig:pce_flowchart} a conceptual scheme of the simulation campaign carried out in this work.

As mentioned, one of the main features of non intrusive PCE is that it can use any deterministic simulation software as a black box input source. We will then present different tests in which PCE has been fed with the output of fluid dynamic simulations based on models characterized by different fidelity levels. In a first test we in fact generated a PCE based on the FOM, and evaluate its performance in a prediction test. The second test consisted in generating a PCE based on the ROM described in the previous section. The latter test allows for an evaluation of how the PCE results are affected when the expansion is based on a surrogate ROM model rather than the FOM one. Given the relatively high number of samples required for the PCE setup, in fact it is interesting to understand if ROM can be used to reduce the computational cost for their generation of the system output at each sample, without a significant loss in terms of accuracy.

\begin{figure}
\centering
\includegraphics[width=0.8\textwidth]{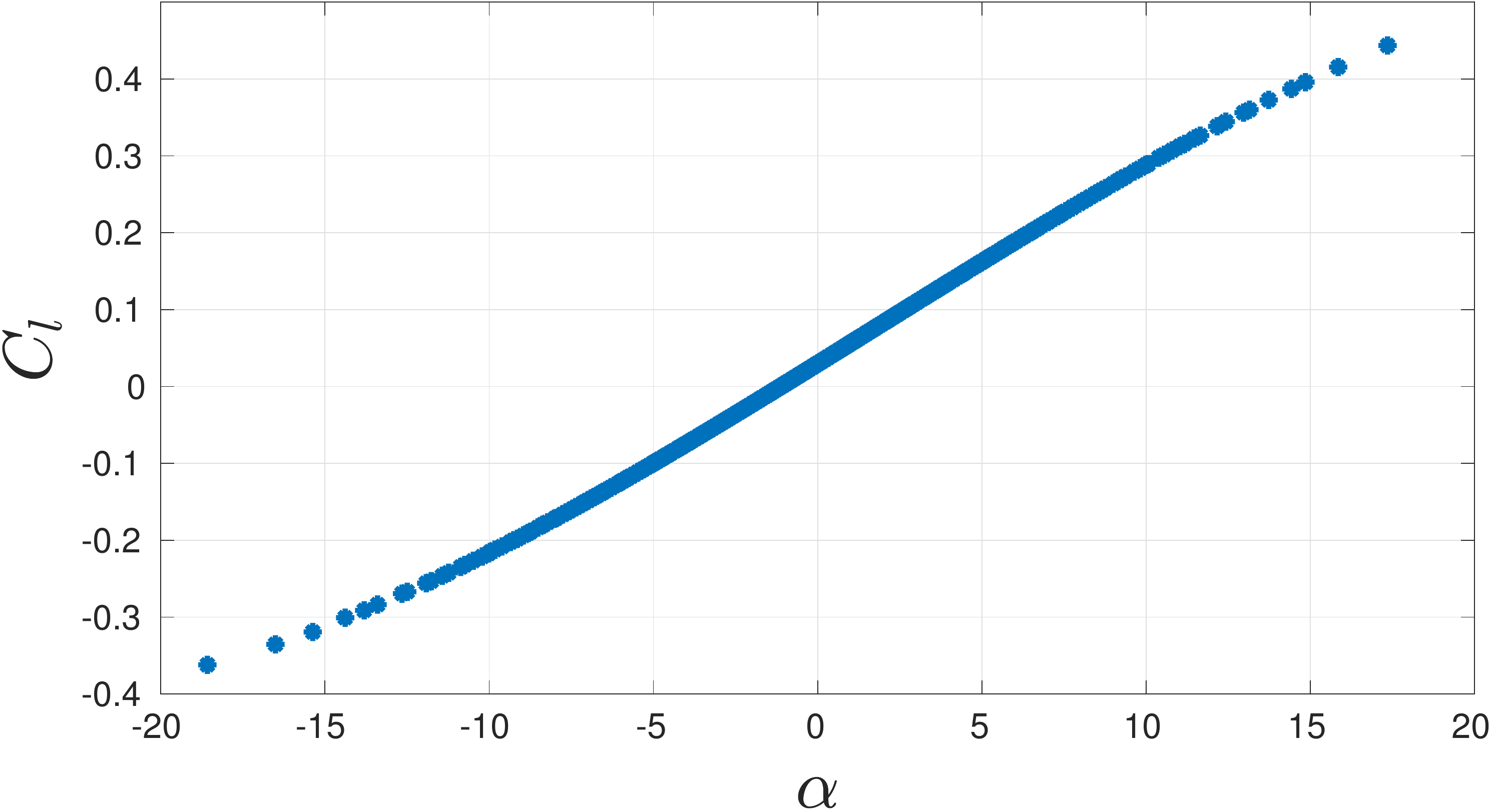}
\caption{The FOM lift coefficient for the UQ case.}\label{fig:LiftHFCase3}
\end{figure}

One of the main assumptions of the non intrusive PCE algorithm implemented is that of operating on Gaussian distributed input parameters. For such reason, the tests in the present section were generated with a set of $300$ sampling points consisting of a single Gaussian distributed bulk. Making use of LHC algorithm, the samples have been randomly generated around angle of attack and velocity magnitude means of $0^{\circ}$ and $100\,\text{m}/\text{s}$ respectively. As for the variances, we prescribed $40^{\circ}$ and $20\,\text{m}/\text{s}$, for angle of attack and velocity magnitude respectively. The FOM lift coefficient curve obtained with the input sampling points described can be seen in \autoref{fig:LiftHFCase3}.\par

\begin{figure}
{
  
    \centering
    \includegraphics[width=0.8\textwidth]{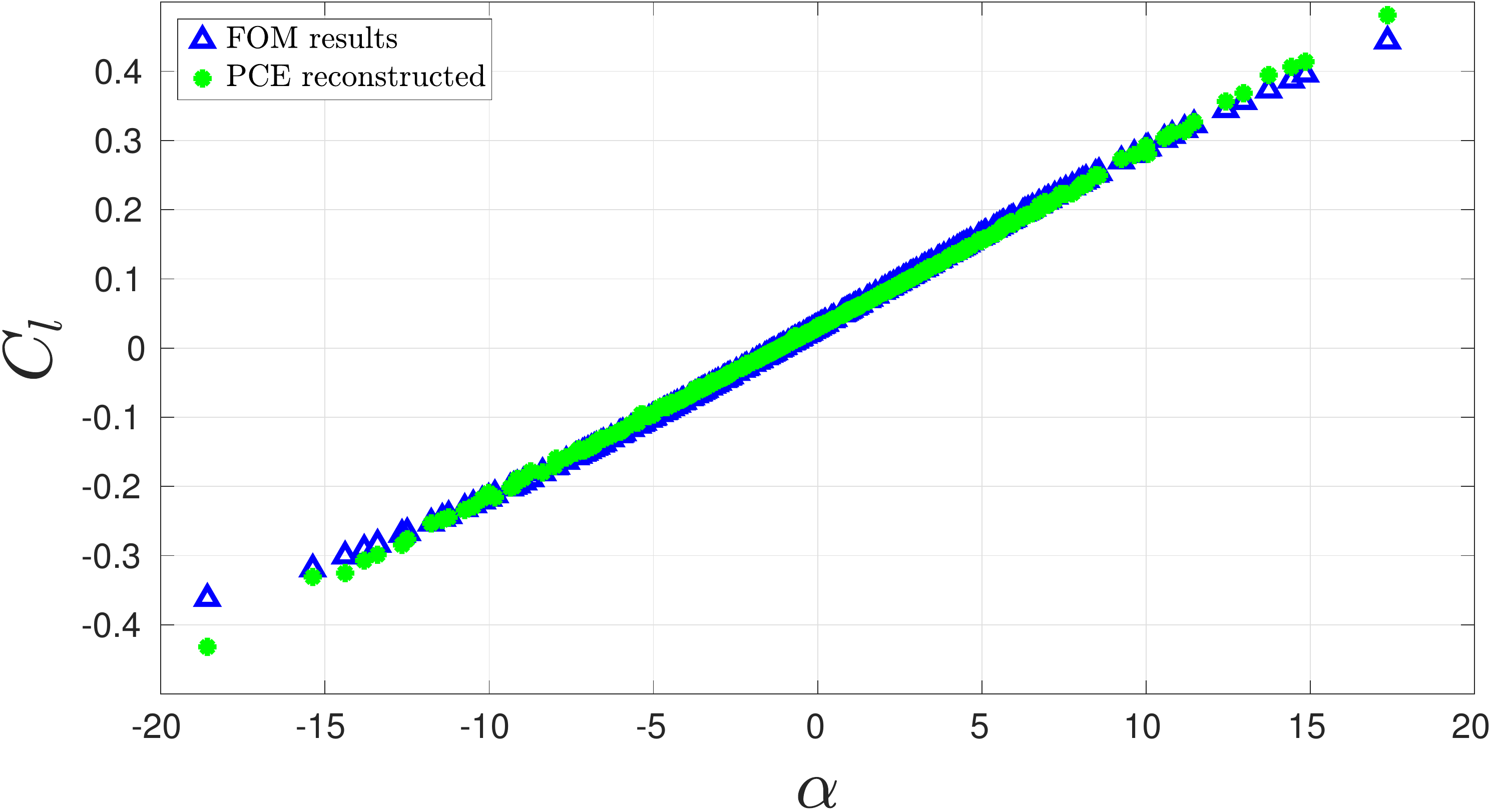}
     \vspace{-0.2cm}\caption{PCE reconstructed lift coefficient versus FOM one, here polynomials of second degree have been used.}\label{fig:LiftComp}
}
\end{figure}

The samples, which have been just mentioned above, will be used as cross validation test for the ROM model developed in the previous subsection. As mentioned earlier, the first test will be to feed PCE with FOM output data and then to conduct a prediction test. PCE is used to predict the value of the lift coefficient for $200$ samples which differ of the $100$ samples used for the PCE coefficients evaluation. Thus the first $100$ samples with their corresponding FOM lift coefficient values were used to build the matrix system \eqref{eq:PCESystem}, which has been solved in the least squares sense. \autoref{fig:LiftComp} displays the $C_l$ values computed with both FOM and PCE in correspondence with all the samples used for check. The overall error in $L^2$ norm between FOM and PCE predications is $5,04\%$.\par  

In the second test we have used ROM data as input for PCE. After using $100$ samples to compute the PCE coefficients, we used the PCE to predict the lift coefficients at $200$ additional samples used for check. We then compared the value of the predicted PCE coefficients in this case to both ROM values and FOM values.\par
\begin{figure}
{
  \begin{minipage}[b]{0.7\linewidth}
    \centering
    \includegraphics[width=\linewidth]{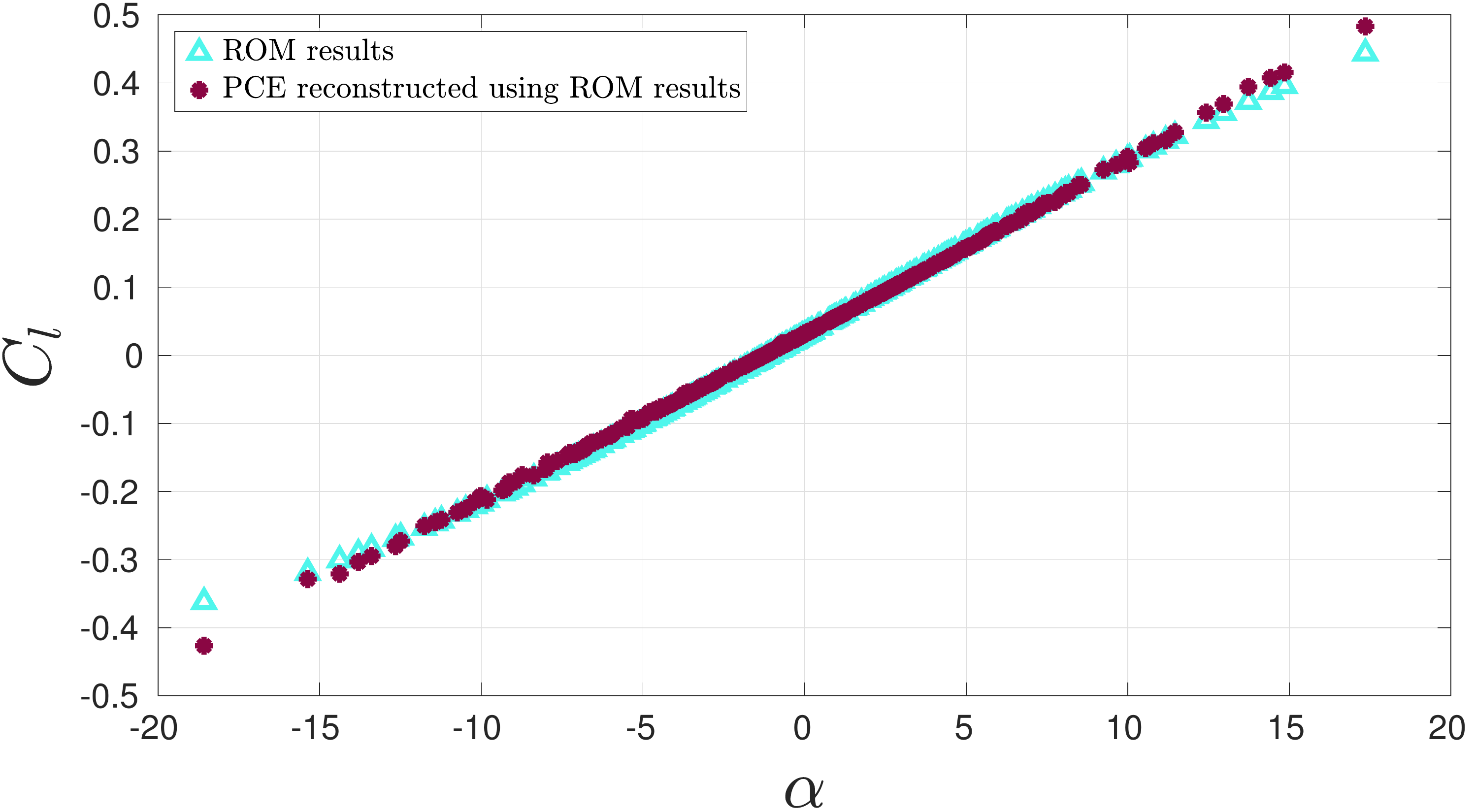}
        \scriptsize(a)

        \end{minipage}
   \begin{minipage}[b]{0.7\linewidth}
    \centering
    \includegraphics[width=\linewidth]{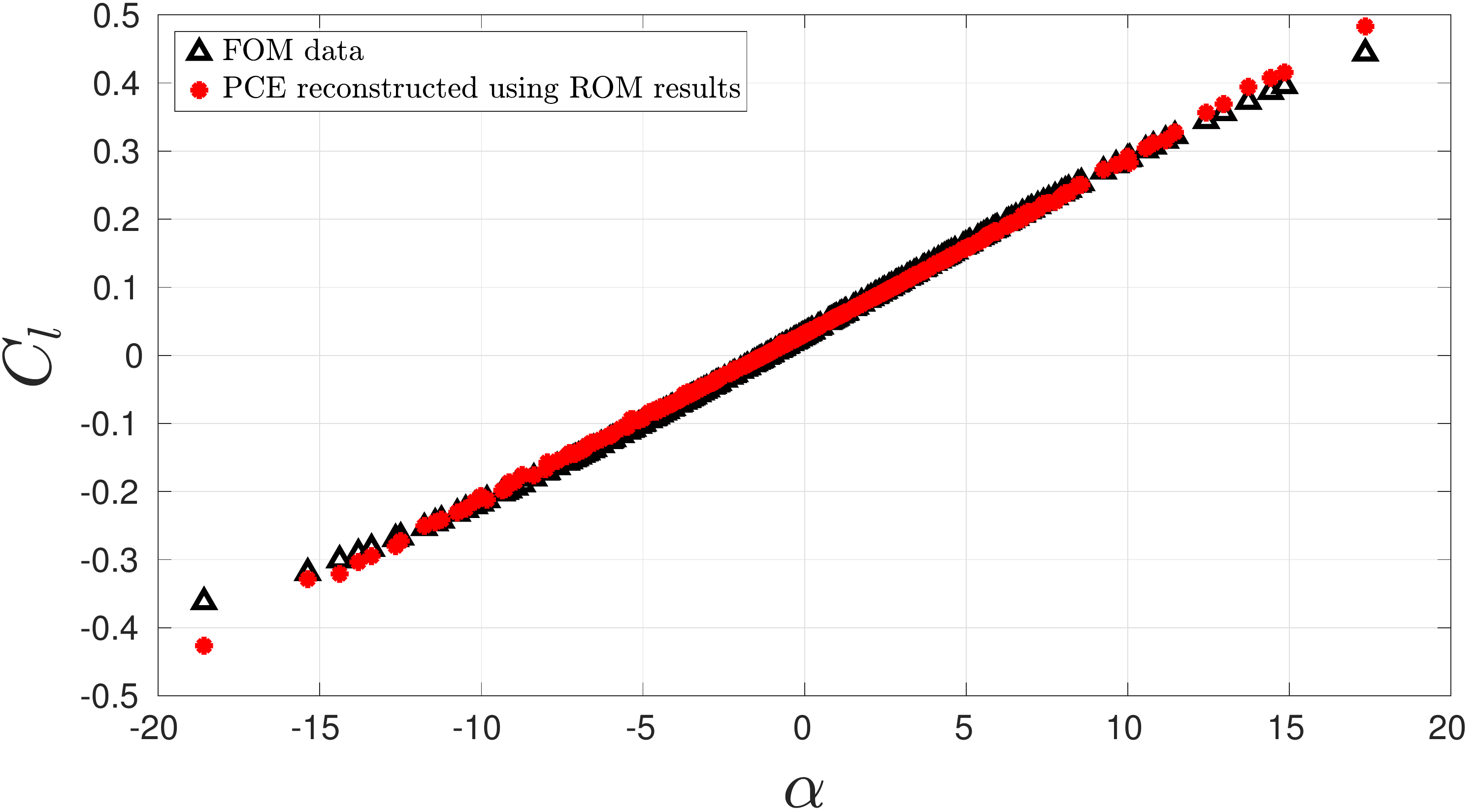}
        \scriptsize(b)
        \end{minipage}     
     \vspace{-0.2cm}\caption{(a) The ROM lift coefficient versus PCE lift coefficient curve when PCE has been applied on ROM output with $30$, $10$ and $10$ modes used for velocity, pressure and supremizer fields respectively. (b) The FOM lift coefficient versus PCE lift coefficient curve when PCE has been applied on ROM output with same number of modes as in (a).}\label{fig:PCEwithPOD}
}
\end{figure}

 
The result of the aforementioned test are reported in \autoref{fig:PCEwithPOD}. The figure includes comparison of the PCE predicted $C_l$ curve with both its ROM and FOM counterparts. The plots show a similar behaviour of the PCE predictions obtained using ROM and FOM output data. By a quantitative standpoint, the PCE predictions present a $4,4\%$ error with respect to the ROM predictions, while the $L^2$ norm of the error with respect to the FOM predictions is $5,14\%$. A summary of the comparisons made is reported in \autoref{table:PCEcomp}.\par 
\begin{table}[htp]
\centering
{
\begin{tabular}{ c | c | c }
First data & Second data & Error \\
\hline  			
FOM  & ROM & $1,53\%$ \\
FOM  & PCE on FOM & $5,04\%$  \\
ROM  & PCE on ROM & $4,4\%$ \\
FOM  & PCE on ROM & $5,14\%$ 
\end{tabular}}\caption{{A comparison between the relative error in $L^2$ norm for the results obtained from ROM and PCE, with PCE being used on both FOM and ROM results. We remark that the number of POD modes used (if apply) are $30, 10$ and $10$ for velocity, pressure and supremizer fields, respectively, for all cases. We underline also that $200$ samples have been used for testing the PCE wherever it is used.}}
\label{table:PCEcomp}
\end{table}
\section{Conclusions and future developments}\label{sec:concl}
In this work, we studied two popular techniques which are used often in the fields of ROM and UQ which are the POD and  PCE, respectively. The study aimed at comparing the accuracy of the two techniques in reconstructing the outputs of interest of viscous fluid dynamic simulations. We have concluded the work with combining the two approaches so as to exploit ROM to speed up the many query problem needed to obtain the PCE coefficients. POD can be a reliable output evaluator for PCE, as the value of relative error PCE had when it was based on ROM results was $5,14\%$ while the error was $5,04\%$ when PCE was based on FOM outputs. The last result speaks positively for POD and makes it a valid tool to be possibly used in the field of uncertainty quantification.\par
The work can be extended in the direction of merging the two approaches in a different way, where one can assume the coefficients in the ROM expansion are not deterministic, but rather dependent on some random variables. The latter assumption can bring UQ into play and one may use techniques such as intrusive/non-intrusive PCE. Our interest is also to extend the proposed methodology, still in the context of reduced order models, to more complex and turbulent flow patterns such as those presented in \cite{HijaziAliStabileBallarinRozza2018,saddam2018,tezzele2018}.  
\section*{Acknowledgments}
We acknowledge the support provided by the European Research Council Consolidator Grant project Advanced Reduced Order Methods with Applications in Computational Fluid Dynamics - GA 681447, H2020-ERC COG 2015 AROMA-CFD (PI: Prof. G. Rozza), MIUR (Italian Ministry of Education, Universities and Research) FARE-X-AROMA-CFD and INdAM-GNCS projects.

\newpage

\bibliographystyle{amsplain_gio}
\bibliography{bibfile_sissa}
\end{document}

%% file: figures/output.pdf_tex
\begingroup%
  \makeatletter%
  \providecommand\color[2][]{%
    \errmessage{(Inkscape) Color is used for the text in Inkscape, but the package 'color.sty' is not loaded}%
    \renewcommand\color[2][]{}%
  }%
  \providecommand\transparent[1]{%
    \errmessage{(Inkscape) Transparency is used (non-zero) for the text in Inkscape, but the package 'transparent.sty' is not loaded}%
    \renewcommand\transparent[1]{}%
  }%
  \providecommand\rotatebox[2]{#2}%
  \ifx\svgwidth\undefined%
    \setlength{\unitlength}{542.50807608bp}%
    \ifx\svgscale\undefined%
      \relax%
    \else%
      \setlength{\unitlength}{\unitlength * \real{\svgscale}}%
    \fi%
  \else%
    \setlength{\unitlength}{\svgwidth}%
  \fi%
  \global\let\svgwidth\undefined%
  \global\let\svgscale\undefined%
  \makeatother%
  \begin{picture}(1,0.60871182)%
    \put(0,0){\includegraphics[width=\unitlength,page=1]{./figures/output.pdf}}%
    \put(0.16438224,0.41574845){\color[rgb]{0,0,0}\makebox(0,0)[lb]{\smash{}}}%
    \put(0,0){\includegraphics[width=\unitlength,page=2]{./figures/output.pdf}}%
    \put(0.27040876,0.35588548){\color[rgb]{0,0,0}\makebox(0,0)[lb]{\smash{$U_\infty$}}}%
    \put(0.43711492,0.46494312){\color[rgb]{0,0,0}\makebox(0,0)[lb]{\smash{$\alpha$}}}%
  \end{picture}%
\endgroup%